\documentclass{amsart}

\usepackage{latexsym}
\usepackage{euscript}
\usepackage{epsfig}

\newenvironment{pf}{\begin{trivlist} \item[] {\it Proof.} \ }{\hfill q.e.d. \end{trivlist} } 
\newenvironment{pf*}[1]{\begin{trivlist} \item[] {\it #1.} \ } {\qed \end{trivlist} }
\newenvironment{ex}{\begin{trivlist} \item[] {\it Example:} }{ \end{trivlist}  }
\newenvironment{contex}{\begin{trivlist} \item[] {\it Example (cont'd):} }{ \end{trivlist}  }

\newtheorem{theorem}{Theorem}[section]
\newtheorem{lemma}[theorem]{Lemma}
\newtheorem{proposition}[theorem]{Proposition}

\newtheorem{theoremvoid}{Theorem}  

\numberwithin{equation}{section}

\renewcommand{\caption}{}

\newcommand{\w}{\wedge}
\newcommand{\om}{\omega}
\newcommand{\alp}{\alpha}
\newcommand{\bet}{\beta}

\newcommand{\lam}{\lambda}

\newcommand{\eps}{\varepsilon}

\newcommand{\calB}{{\mathcal B}}
\newcommand{\calC}{{\mathcal C}}
\newcommand{\calI}{{\mathcal I}}
\newcommand{\calJ}{{\mathcal J}}
\newcommand{\calK}{{\mathcal K}}
\newcommand{\calL}{{\mathcal L}}
\newcommand{\calF}{{\mathcal F}}

\newcommand{\hook}{\, {\rule[0in]{2mm}{0.25mm}\rule{0.25mm}{2mm}} \, }
\newcommand{\bb}{\mathbb}
\newcommand{\zbar}{\bar{z}}
\def\B/{B\"acklund}

\begin{document}

\title[Homogeneous B\"acklund transformations]
{Homogeneous B\"acklund transformations of hyperbolic Monge-Amp\`ere systems}

\author {Jeanne Nielsen Clelland}
\address{Department of Mathematics, 395 UCB, University of 
Colorado\\ 
Boulder, CO 80309-0395}
\email{Jeanne.Clelland@colorado.edu}

\thanks{This research was supported in part by NSF grant DMS-9627403.}

\begin{abstract}
A B\"acklund transformation between two hyperbolic Monge-Am\-p\`ere systems 
may be described as a certain type of exterior 
differential system on a 6-dimensional manifold $\calB$.  The transformation 
is {\em homogeneous} if the group of symmetries of the system acts 
transitively on $\calB$.
We give a complete classification of homogeneous B\"acklund 
transformations between hyperbolic Monge-Amp\`ere systems.  
\end{abstract}

\maketitle

\section{Introduction}

In this paper we will study B\"acklund transformations between two 
hyperbolic Monge-Amp\`ere equations.  A {\em Monge-Amp\`ere equation} 
is a partial differential equation of the form
\[ A (z_{xx}z_{yy} - z_{xy}^2) + B z_{xx} + 2C z_{xy} + D z_{yy} + 
E = 0 \]
where the coefficients $A,B,C,D,E$ are functions of the variables 
$x,y,z,z_x, z_y$.  The equation is {\em hyperbolic} if it has 
distinct, real characteristics at each point, i.e., if $AE - BD + 
C^2 > 0$.

There are many definitions of
B\"acklund transformations given in the literature.  Rather than 
attempting
to give an all-encompassing definition, we will use B\"acklund's 
original notion.  In \cite{B80} he posed the following general problem:
Let $M^5 = \bar{M}^5 = \bb{R}^5$, with 
coordinates $(x,y,z,p,q)$ on $M$ and 
$(\bar{x},\bar{y},\bar{z},\bar{p},\bar{q})$ on $\bar{M}$.  Given four 
equations
\begin{equation} 
F_i(x,y,z,p,q,\bar{x},\bar{y},\bar{z},\bar{p},\bar{q}) = 0, \ \ \ 
i=1,\ldots, 4, \label{foureqs}
\end{equation}
find parametrized surfaces 
$X: U \to M, \ \bar{X}: U \to \bar{M}$
and a one-to-one correspondence between them (which may be 
established by using the same parameter domain $U$ for both surfaces)
such that 
the coordinate functions $x(u,v), \ldots, \bar{q}(u,v)$  
of the two surfaces satisfy the conditions
\begin{gather*}
F_i(x(u,v), \ldots, \bar{q}(u,v)) = 0, \ \ \ i=1, \ldots, 4 \\
dz = p\,dx + q\,dy  \\
d\bar{z} = \bar{p}\,d\bar{x} + \bar{q}\,d\bar{y}. 
\end{gather*}
(The last two equations imply that the coordinates $p,q,\bar{p}, 
\bar{q}$ should be regarded as the partial derivatives $z_x, z_y, 
\bar{z}_{\bar{x}}, \bar{z}_{\bar{y}}$, respectively.)

B\"acklund's approach to this problem was to assume that $X$ is a 
graph of the form
\[ (x,y,z,p,q) = (x,y, z(x,y), z_x(x,y), z_y(x,y)) \]
for some known function $z(x,y)$.
Two of the equations \eqref{foureqs} can be solved for the variables $x$ and 
$y$, and substituting these expressions into the remaining two 
equations yields equations of the form
\begin{gather}
f(\bar{x},\bar{y},\bar{z},\bar{p},\bar{q}) = 0 \label{back} \\
g(\bar{x},\bar{y},\bar{z},\bar{p},\bar{q}) = 0. \notag
\end{gather}
This may be regarded as an overdetermined first-order PDE system for the 
function $\bar{z}(\bar{x}, 
\bar{y})$, and the compatibility conditions for this system take the 
form of partial differential equations that must be satisfied by the 
function $z(x,y)$.  
If $z(x,y)$ satisfies these conditions, then the system \eqref{back} 
has a 1-parameter family of solutions 
$\bar{z}(\bar{x}, \bar{y})$ which can be found by solving ordinary 
differential equations.  In this case, the four equations 
\eqref{foureqs} may 
be regarded as a transformation of the surface $z=z(x,y)$ into the
surface $\bar{z} = \bar{z}(\bar{x}, \bar{y})$.

\begin{ex}
The classical B\"acklund transformation for the 
sine-Gordon equation
\begin{equation}
 z_{xy} = \tfrac{1}{2}\sin (2z) \label{SGE}
\end{equation}
is usually described by the two equations
\begin{align}
z_x + \zbar_x & = \lam\, \sin(z - \zbar) \label{SGEback} \\
z_y - \zbar_y & = \frac{1}{\lam}\, \sin(z + \zbar) \notag
\end{align}
where $\lam$ is a nonzero constant.
In B\"acklund's notation we would write equations \eqref{SGEback} as
\begin{align*}
p + \bar{p} & = \lam\, \sin(z - \zbar) \label{SGEback} \\
q - \bar{q} & = \frac{1}{\lam}\, \sin(z + \zbar),
\end{align*}
together with the two additional equations $\bar{x} = x,\, \bar{y} = y$.
It is straightforward to show that if 
$z(x,y), \zbar(x,y)$ satisfy equations \eqref{SGEback}, then both 
must be solutions of the sine-Gordon equation \eqref{SGE}.  Moreover, 
if $z(x,y)$ is any known solution of \eqref{SGE}, then \eqref{SGEback} is a
compatible, overdetermined system for the unknown function 
$\zbar(x,y)$, whose solution depends only on solving ordinary 
differential equations.  For instance, taking $z(x,y) = 0$ yields
\[ \zbar(x,y) = \tan^{-1}(e^{-(\lam x + \tfrac{1}{\lam}y + c)}). \]
These are the {\em 1-soliton} solutions of the sine-Gordon equation.
\end{ex}

The problem that arises in B\"acklund's approach is that for a generic 
choice of equations \eqref{foureqs}, the compatibility conditions for 
\eqref{back} cannot be written as separate PDEs for the functions $z$ 
and $\bar{z}$; rather they involve $z$ and $\bar{z}$ together, and 
so equations \eqref{foureqs} do not give a transformation of the desired form.  
This raises the question: for what sets of equations \eqref{foureqs} 
can the compatibility conditions for \eqref{back} be written as separate 
PDEs for $z$ and $\bar{z}$?  If this is the case, then the system 
\eqref{foureqs} is called a 
{\em B\"acklund transformation} between the two PDEs.  Thus the 
question raised above becomes: what 
PDEs (or pairs of PDEs) have B\"acklund transformations?
This is an open problem which has attracted much attention over the 
past century.  Its difficulty is attested to by the
extensive work that analysts such as Goursat \cite{G21,G22} put into
investigating
special cases, such as equations of the form
$$z_{xy}=\rho z_x z_y +a z_x + b z_y + c$$ where $a,b,c,\rho$ are
functions of $x,y$ and $z$.  
More recently, McLaughlin and Scott \cite{MS73} classified
auto-B\"acklund transformations (i.e., transformations for which $z$ 
and $\bar{z}$ satisfy the same PDE) of equations of the form
\[ z_{xy} + a z_x + b z_y = F(z) \]
where $a$ and $b$ are constants, and Byrnes \cite{B76} generalized 
this work by allowing $F$ to depend on $x$ and $y$ as well as $z$.  
Zvyagin \cite{Z81,Z86}, following Goursat's 
approach, has 
studied a certain type of B\"acklund transformation which he calls 
{\em harmonic}; he has also given a classification \cite{Z91} of 
B\"acklund transformations of the wave equation $z_{xy}=0$, 
although the descriptions of the systems on his list are somewhat 
unsatisfying and his paper contains no proof.  These references 
represent only a small sample of the work that has been done on this 
problem; it would be impossible to give a complete list.

Although 
Goursat's foundational work appears to be highly dependent
on working in coordinates, he was the first to focus on the geometric
structures underlying B\"acklund transformations.  This approach has 
since proven quite fruitful, and
these structures are best decribed in terms of exterior 
differential systems.  

An {\em exterior differential system} on a manifold $M$ is a differentially
closed ideal ${\calI}$ in the algebra of differential
forms on $M$.  Any system of partial differential equations can be 
formulated as an exterior differential
system $\calI$, and solutions of the PDE system correspond to {\em integral
manifolds} of $\calI$, i.e., submanifolds $N \subset M$ which
satisfy the condition that all the forms in $\calI$ vanish when restricted
to $N$.  
A {\em Monge-Amp\`ere system} $\calI$ is an exterior differential
system on a 5-dimensional manifold $M$ that is locally generated by a
contact form $\theta$ (i.e., a 1-form $\theta$ with the property that 
$\theta \w d\theta \w d\theta \neq 0$), the 2-form $\Theta = d\theta$, 
and another 2-form $\Psi$.
A Monge-Amp\`ere system $\calI$ is {\em hyperbolic} if
the quadratic equation
\[ (\lambda\, \Theta + \mu\, \Psi) \w (\lambda\, \Theta + \mu\, \Psi) 
\equiv 0 \mod{\theta} \]
has distinct, real roots.  This condition agrees with the traditional 
definition of hyperbolicity, and it implies that there are two 
independent linear combinations $\lambda\, \Theta + \mu\, \Psi$ which are 
{\em decomposable} 2-forms (i.e., 2-forms which can be written as 
$\om^1 \w \om^2$ for some 1-forms $\om^1, \om^2$) modulo 
$\theta$.  (See \cite{BGH95} for a 
discussion of hyperbolic exterior differential systems.)

\begin{contex}
The sine-Gordon equation \eqref{SGE} may be described 
as a hyperbolic Monge-Amp\`ere system on $\mathbb{R}^5$ (with 
coordinates $(x,y,z,p,q)$)
generated by the forms
\begin{gather*}
 \theta = dz - p\,dx - q\, dy \notag \\  
 \Theta = -dp \w dx - dq \w dy \\
 \Psi = [dp - \tfrac{1}{2} \sin(2z)\, dy] \w dx . \notag
\end{gather*}
Note that $\Psi$ is decomposable; the other decomposable linear 
combination of $\Psi$ and $\Theta$ is $-(\Psi + \Theta) = [dq - \tfrac{1}{2} 
\sin(2z)\, dx] \w dy$.
Two-dimensional integral manifolds of this system that satisfy the 
independence condition
$dx \w dy \neq 0$ are naturally in one-to-one correspondence with solutions of
\eqref{SGE}.
\end{contex}

B\"acklund's original notion may 
be expressed in this context as follows.
Suppose that $(M_1, \calI_1)$ and $(M_2, \calI_2)$ are hyperbolic
Monge-Amp\`ere systems, with
\[ \calI_1 = \{ \theta_1, \Theta_1, \Psi_1 \} \]
\[ \calI_2 = \{ \theta_2, \Theta_2, \Psi_2 \}.\]
A {\em B\"acklund transformation} between $(M_1, \calI_1)$ and $(M_2, 
\calI_2)$ is a 6-dimensional submanifold $\calB \subset M_1 \times M_2$
which has the following properties:
\begin{enumerate}
\item{The natural projections $\pi_1:\calB \rightarrow M_1$ and $\pi_2: \calB
\rightarrow M_2$ are submersions.

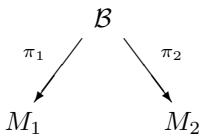
\setlength{\unitlength}{2pt}
\begin{center}
\begin{picture}(40,30)(0,0)
\put(5,5){\makebox(0,0){$M_1$}}
\put(35,5){\makebox(0,0){$M_2$}}
\put(20,25){\makebox(0,0){$\calB$}}
\put(16,21){\vector(-3,-4){9}}
\put(24,21){\vector(3,-4){9}}
\put(7,18){\makebox(0,0){$\scriptstyle{\pi_1}$}}
\put(33,18){\makebox(0,0){$\scriptstyle{\pi_2}$}}
\end{picture}
\end{center}  
}
\item{The pullbacks to $\calB$ of the forms $\Theta_1, \Theta_2, 
\Psi_1, \Psi_2$ satisfy the condition that
\begin{equation*}
 \{ \Psi_1, \Psi_2 \} \equiv \{ \Theta_1, \Theta_2 \}  
   \mod{\{\theta_1, \theta_2\}}. 
\end{equation*}
Since $\Theta_1, \Psi_1$ are linearly independent forms (as are
$\Theta_2, \Psi_2$), this condition implies that
\begin{equation*}
 \{ \Theta_1, \Psi_1 \}  \equiv \{ \Theta_2, \Psi_2  \} 
   \mod{\{\theta_1, \theta_2\}}. 
\end{equation*}
This second equation is really the desired
property; the first equation ensures that, in 
addition, the forms $\Theta_1, \Theta_2$ are linearly independent.  }
\end{enumerate}

That this definition captures the desired behavior may be seen as 
follows: suppose that $N
\hookrightarrow M_1$ is a 2-dimensional integral manifold of $\calI_1$.  
The inverse image
$\pi_1^{-1}(N)$ is a 3-dimensional submanifold of $\calB$.  Now consider the
restriction of $\pi_2^* (\calI_2)$ to $\pi_1^{-1}(N)$.  By 
Property (2) above, the restriction of $\pi_2^* (\calI_2)$ is a
Frobenius system (i.e., an exterior differential system which is 
generated {\em algebraically} by its 1-forms) on $\pi_1^{-1}(N)$.   
By the Frobenius Theorem, $\pi_1^{-1}(N)$ is foliated by 2-dimensional
integral manifolds of $\pi_2^* (\calI_2)$, each of which 
projects to an integral manifold of $(M_2, \calI_2)$; moreover, these 
integral manifolds
can be constructed by solving ODEs.

From the point of view of B\"acklund's original problem, any 2-dimensional 
integral manifold 
$\tilde{S} \subset \calB$ of the ideal $\calJ = \{\theta_1, \theta_2, \Theta_1, 
\Theta_2 \}$ projects to surfaces $S_1 \subset M_1, \ S_2 \subset 
M_2$ which 
are integral manifolds of $\calI_1, \calI_2$ respectively.  The 
condition that $\tilde{S}$ be an integral manifold of $\calJ$ is 
exactly the requirement that the compatibility conditions for the
equations \eqref{back} be satisfied.

Our primary tool for classifying such structures will be Cartan's method of 
equivalence; this is a method for computing local invariants
of exterior differential systems and deciding when two systems
are equivalent under some natural class of diffeomorphisms. In principle, 
it should be possible to completely 
classify all B\"acklund transformations of hyperbolic Monge-Amp\`ere 
systems using this method.  Unfortunately, in practice it is rarely 
possible to carry out this process in full generality.  In this paper 
we will perform the somewhat simpler task of classifying the {\em 
homogeneous} B\"acklund transformations, i.e., those transformations for 
which the group of symmetries of the structure
$(\calB,\ \calI_1,\  \calI_2)$ acts transitively on $\calB$.
The main result is the following theorem.

\begin{theoremvoid}[cf. Theorem \ref{bigthm}]
Let $\calB \subset M_1 \times M_2$ be a homogeneous B\"acklund 
transformation.  Then $\calB$ is locally contact 
equivalent to one of the following:
\begin{enumerate}
\item{A B\"acklund transformation between solutions of the wave 
equation $z_{xy} = 0$}
\item{A holonomic B\"acklund transformation of the form described in Theorem 
\ref{case3athm}}
\item{The classical B\"acklund transformation between the wave 
equation $z_{xy} = 0$ and Liouville's equation $z_{xy} = e^z$}
\item{A B\"acklund transformation between surfaces of constant 
negative Gauss curvature in $\bb{E}^3$}
\item{A B\"acklund transformation between surfaces of constant 
Gauss curvature $0 < K < 1$ in $S^3$}
\item{A B\"acklund transformation between surfaces of constant 
Gauss curvature $-\infty < K < -1$ in $\bb{H}^3$}
\item{A B\"acklund transformation between spacelike surfaces of 
constant positive Gauss curvature in $\bb{E}^{2,1}$}
\item{A B\"acklund transformation between timelike surfaces of 
constant positive \linebreak Gauss curvature, or equivalently, constant nonzero 
mean curvature, in $\bb{E}^{2,1}$}
\item{A B\"acklund transformation between timelike minimal surfaces 
in $\bb{E}^{2,1}$}
\item{A B\"acklund transformation between spacelike surfaces of 
constant Gauss curvature $1 < K < \infty$ in $S^{2,1}$}
\item{A B\"acklund transformation between timelike surfaces of 
constant Gauss curvature $1 < K < \infty$, or equivalently, constant 
mean curvature $H \in \bb{R}$, in $S^{2,1}$}
\item{A B\"acklund transformation between spacelike surfaces of 
constant Gauss curvature $-1 < K < \infty, \ K \neq 0$ in $\bb{H}^{2,1}$}
\item{A B\"acklund transformation between timelike surfaces of 
constant Gauss curvature $-1 < K < \infty, \ K \neq 0$, or equivalently, 
constant mean curvature $|H| > 1$, in $\bb{H}^{2,1}$}
\item{A B\"acklund transformation between timelike surfaces of 
constant mean curvature $|H| \leq 1$ in $\bb{H}^{2,1}$.}
\item{A B\"acklund transformation between certain surfaces in a 
5-dimensional quotient space of $SO^{\ast}(4)$.}
\end{enumerate}
\end{theoremvoid}

Throughout this paper we will work locally.  Statements such as 
``assume that $C \neq 0$'' should be interpreted as ``assume that $C$ 
is not identically zero and restrict to the open set where $C \neq 0$''.

\section{The equivalence problem}

Suppose that $\calB$ is a B\"acklund transformation between two 
hyperbolic Monge-Amp\`ere systems $(M_1, \calI_1)$ and $(M_2, 
\calI_2)$.  Let $\calJ$ be the ideal on $\calB$ generated by the 
pullbacks of $\calI_1$ and $\calI_2$; according to our definition of a 
B\"acklund transformation, $\calJ$ is generated algebraically by the forms
$\{\theta_1, \theta_2, \Theta_1, \Theta_2 \}. $ 

Since $\calI_1$ and $\calI_2$ are hyperbolic, locally there exist 1-forms 
$\om^1, \om^2, 
\om^3, \om^4$ on $\calB$ such that $\{\theta_1,\, \theta_2,\, \om^1,\, 
\om^2,\, \om^3,\, \om^4\}$ is a coframing of $\calB$ (i.e., a basis for 
the space of 1-forms on $\calB$)  and 
\[ \calJ = \{\theta_1,\, \theta_2,\, \om^1 \w \om^2,\, \om^3 \w \om^4\}. 
\]
(It is important to note that $\theta_1$ and $\theta_2$ are each 
{\em separately} determined up to a scalar multiple, since $\theta_i$ 
determines the contact structure on $M_i$.)
Any such coframing has the property that
\begin{align*}
d\theta_1 & \equiv A_1\, \om^1 \w \om^2 + A_2\,\om^3 \w \om^4  
\mod{\{\theta_1,\, \theta_2\}}\\
d\theta_2 & \equiv A_3\, \om^1 \w \om^2 + A_4\, \om^3 \w \om^4  
\mod{\{\theta_1,\, \theta_2\}}
\end{align*}
for some nonvanishing functions $A_1, A_2, A_3, A_4$.  Since $d\theta_1,\, 
d\theta_2$ are required to be linearly independent 2-forms at each 
point of $\calB$, we must have $A_1 A_4 - A_2 A_3 \neq 0$.

By rescaling the $\om^i$ and adding multiples of $\theta_1$ and 
$\theta_2$ to the $\om^i$ if necessary, we can arrange that
\begin{align}
d\theta_1 & \equiv A_1\, \om^1 \w \om^2 + \om^3 \w \om^4  
\mod{\theta_1} \label{dtheta}\\
d\theta_2 & \equiv \om^1 \w \om^2 + A_2\, \om^3 \w \om^4  
\mod{\theta_2} \notag
\end{align}
for some nonvanishing functions $A_1,\, A_2$ on $\calB$ with $A_1 
A_2 \neq 1$.
This coframing is not unique; any other such coframing 
$\{\tilde{\theta}_1,\, \tilde{\theta}_2,\, \tilde{\om}^1,\, 
\tilde{\om}^2,\, 
\tilde{\om}^3,\, \tilde{\om}^4\}$ has the form
\begin{equation}
 \begin{bmatrix} \tilde{\theta}_1\\[0.1in] \tilde{\theta}_2\\[0.1in] 
\tilde{\om}^1\\[0.1in] \tilde{\om}^2\\[0.1in] \tilde{\om}^3\\[0.1in] 
\tilde{\om}^4 
\end{bmatrix} =   \begin{bmatrix} b_{11}b_{22} - b_{12}b_{21} & 0 & 
0 & 0 & 0 & 0 \\[0.1in] 0 &  a_{11}a_{22} - a_{12}a_{21} & 0 & 0 & 0 & 0 
\\[0.1in] 0 & 0 & a_{11} & a_{12} & 0 & 0 \\[0.1in] 0 & 0 & a_{21} & 
a_{22} & 0 & 0 
\\[0.1in] 0 & 0 & 0 & 0 & b_{11}&  b_{12} \\[0.1in] 0 & 0 & 0 & 
0 & b_{21} & b_{22} 
\end{bmatrix}^{-1} \begin{bmatrix} \theta_1\\[0.1in] \theta_2\\[0.1in] 
\om^1\\[0.1in] 
\om^2\\[0.1in] \om^3\\[0.1in] \om^4 \end{bmatrix}  \label{G0freedom}
\end{equation}
where 
$ b_{11}b_{22} - b_{12}b_{21} \neq 0, \ 
a_{11}a_{22} - a_{12}a_{21} \neq 0.$  (The inverse is included for 
greater ease of computation in what follows.)  A coframing satisfying 
\eqref{dtheta} is called {\em 0-adapted}, and the group $G_0$ of matrices 
of the 
above form is called the {\em structure group} of the equivalence 
problem.  (In fact, the most general choice of structure group would 
include a discrete component interchanging the 
distributions $\{\om^1, \om^2\}$ and $\{\om^3, \om^4\}$.  However, 
this freedom does not contribute anything crucial to the structure 
group, and it is easier to work with a connected group.)

Now consider the exterior derivatives of the $\om^i$.  Because 
$\theta_1$ is well-defined (up to scalar multiples) on 
$M_1$, its Cartan system 
$\calC = \{\theta_1, \om^1, \om^2, \om^3, \om^4\}$ is well-defined on 
$M_1$.  (The {\em Cartan system} of a 1-form $\theta$ may be thought of as the 
span of a minimal set of 1-forms required to express $\theta$ and 
$d\theta$.  It is always a Frobenius system; see \cite{BCG3} 
for details.)  
  In fact, $M_1$ is (locally) the quotient of $\calB$ by the leaves of the 
foliation defined by $\calC$.  
Let $\{\frac{\partial}{\partial \theta_1}, \frac{\partial}{\partial 
\theta_2}, \frac{\partial}{\partial \om^1}, \frac{\partial}{\partial 
\om^2}, \frac{\partial}{\partial \om^3}, \frac{\partial}{\partial 
\om^4}\}$ denote the basis for the tangent space of $\calB$ which is 
dual to the coframing $\{\theta_1,\, \theta_2,\, \om^1,\, \om^2,\, 
\om^3,\, \om^4\}$.
The ideal 
\[ \calI_1 = \{\theta_1, \om^1 \w \om^2, \om^3 \w \om^4 \} \]
is well-defined on $M_1$, as are its characteristic systems 
\[ \calK_{11} = \{\theta_1, \om^1 \w \om^2\}, \qquad \calK_{12} = \{\theta_1, 
\om^3 \w \om^4\}. \]
Therefore the Lie derivative $\calL_{\frac{\partial}{\partial 
\theta_2}} (\om^1 \w \om^2)$ must satisfy
\begin{align}
0 &\equiv \calL_{\frac{\partial}{\partial \theta_2}} (\om^1 \w \om^2) 
\mod{\{\theta_1, \om^1 \w \om^2\}} \notag \\
&\equiv \frac{\partial}{\partial \theta_2} \hook (d\om^1 \w \om^2 - 
\om^1 \w d\om^2)\mod{\{\theta_1, \om^1 \w \om^2\}}. \label{Lieeq1}
\end{align}
Reducing equation \eqref{Lieeq1} modulo $\om^1$ yields
\[ \left(\frac{\partial}{\partial \theta_2} \hook d\om^1 \right) \w 
\om^2 \equiv 0 \mod{\{\theta_1, \om^1\}}, \]
and therefore
\[ \frac{\partial}{\partial \theta_2} \hook d\om^1 \equiv 0 
\mod{\{\theta_1, \om^1, \om^2\}}. \]
Consequently, $d\om^1$ cannot contain any terms involving the 
2-forms $\theta_2 \w \om^3$ or $\theta_2 \w \om^4$.  Similarly, reducing 
equation \eqref{Lieeq1} modulo $\om^2$ shows that $d\om^2$ cannot contain 
any terms involving the 2-forms $\theta_2 \w \om^3$ or $\theta_2 \w \om^4$.
An analogous argument using the equation 
\[ \calL_{\frac{\partial}{\partial \theta_2}} (\om^3 \w \om^4) \equiv 
0 \mod{\{\theta_1, \om^3 \w \om^4\}} \]
shows that $d\om^3$ and $d\om^4$ cannot contain any terms involving the 
2-forms $\theta_2 \w \om^1$ or $\theta_2 \w \om^2$.

This argument can be repeated for the characteristic systems 
\[ \calK_{21} = \{\theta_2, \om^1 \w \om^2\}, \qquad \calK_{22} = 
\{\theta_2, \om^3 \w \om^4\} \]  
of $(M_2, \calI_2)$; this shows that $d\om^1$ and $d\om^2$ 
cannot contain any terms involving the 2-forms  $\theta_1 \w \om^3$ or 
$\theta_1 \w \om^4$, and $d\om^3$ and $d\om^4$ cannot contain any terms 
involving the 2-forms 
$\theta_1 \w \om^1$ or $\theta_1 \w \om^2$.  It follows that
that 
\begin{align*}
& \left. \begin{array}{l}
d\om^1 \equiv B_1\, \theta_1 \w \theta_2 + C_1\, \om^3 \w \om^4 
\\[0.1in]
d\om^2 \equiv B_2\, \theta_1 \w \theta_2 + C_2\, \om^3 \w \om^4
\end{array} \right\} \mod{\{\om^1, \om^2\}} \\[0.1in]
& \left. \begin{array}{l}
d\om^3 \equiv B_3\, \theta_1 \w \theta_2 + C_3\, \om^1 \w \om^2 
\\[0.1in]
d\om^4 \equiv B_4\, \theta_1 \w \theta_2 + C_4\, \om^1 \w \om^2
\end{array} \right\} \mod{\{\om^3, \om^4\}}
\end{align*}
for some functions $B_i, C_i$ on $\calB$.

These equations, 
taken together with equations \eqref{dtheta}, form the {\em structure 
equations}
\begin{equation}
\begin{bmatrix} d\theta_1\\[0.1in] d\theta_2\\[0.1in] d\om^1\\[0.1in] 
d\om^2\\[0.1in] d\om^3\\[0.1in] d\om^4
\end{bmatrix} = 
-\begin{bmatrix}\bet_1 + \bet_4 & 0 & 0 & 0 & 0 & 0 \\[0.1in] 0 & \alp_1 + 
\alp_4 & 0 & 0 & 0 & 0 \\[0.1in] 0 & 0 & \alp_1 & \alp_2 & 0 & 0 
\\[0.1in] 0 & 0 & 
\alp_3 & \alp_4 & 0 & 0 \\[0.1in] 0 & 0 & 0 & 0 & \bet_1 & \bet_2 
\\[0.1in] 0 & 0 & 0 & 
0 & \bet_3 & \bet_4 \end{bmatrix} \w
\begin{bmatrix} \theta_1\\[0.1in] \theta_2\\[0.1in] \om^1\\[0.1in] 
\om^2\\[0.1in] \om^3\\[0.1in] \om^4 
\end{bmatrix} + 
\begin{bmatrix} \Theta_1\\[0.1in] \Theta_2\\[0.1in] \Omega^1\\[0.1in] 
\Omega^2\\[0.1in] \Omega^3\\[0.1in] \Omega^4 
\end{bmatrix} \label{G0struct}
\end{equation}

where the $\alp_i, \bet_i$ are 1-forms on $\calB$ and 
\begin{align*}
\Theta_1 & = \gamma \w \theta_1 + A_1\, \om^1 \w \om^2 + \om^3 \w \om^4  
\\
\Theta_2 & = \delta \w \theta_2 + \om^1 \w \om^2 + A_2\, \om^3 \w 
\om^4 \\
\Omega^1 & = B_1\, \theta_1 \w \theta_2 + C_1\, \om^3 \w \om^4  \\
\Omega^2 & = B_2\, \theta_1 \w \theta_2 + C_2\, \om^3 \w \om^4  \\
\Omega^3 & = B_3\, \theta_1 \w \theta_2 + C_3\, \om^1 \w \om^2  \\
\Omega^4 & = B_4\, \theta_1 \w \theta_2 + C_4\, \om^1 \w \om^2
\end{align*}
for some 1-forms $\gamma, \delta$ on $\calB$.
These equations are chosen so 
that the matrix in \eqref{G0struct} takes values in the Lie algebra 
$\mathfrak{g}_0$
of $G_0$; this is in accordance with the method of equivalence.  (See 
\cite{G89} for details.)
The functional coefficients of the terms appearing in $\Theta_i, 
\Omega^i$ are called {\em torsion} terms.

 We can modify the $\alp_i, 
\bet_i$ if necessary to arrange that
\begin{align*}
\gamma & = E_1\, \theta_2 + F_1\, \om^1 + F_2\, \om^2 \\
\delta & = E_2\, \theta_1 + F_3\, \om^3 + F_4\, \om^4
\end{align*}  
for some functions $E_i, F_i$ on $\calB$.
The forms $\alp_i, \bet_i$ are still not uniquely determined; they are 
determined only up to transformations of the form
\begin{alignat}{2}
\alp_1 & \mapsto \alp_1 + r_1\, \om^1 + r_2\, \om^2 \qquad & \bet_1 & 
\mapsto \bet_1 + s_1\, \om^3 + s_2\, \om^4 \notag \\
\alp_2 & \mapsto \alp_2 + r_2\, \om^1 + r_3\, \om^2 \qquad & \bet_2 & 
\mapsto \bet_2 + s_2\, \om^3 + s_3\, \om^4 \label{G0connfreedom} \\
\alp_3 & \mapsto \alp_3 + r_4\, \om^1 - r_1\, \om^2 \qquad & \bet_3 & 
\mapsto \bet_3 + s_4\, \om^3 - s_1\, \om^4 \notag \\
\alp_4 & \mapsto \alp_4 - r_1\, \om^1 - r_2\, \om^2 \qquad & \bet_4 & 
\mapsto \bet_4 - s_1\, \om^3 - s_2\, \om^4.\notag 
\end{alignat}

Differentiating the structure equations yields
\begin{alignat*}{2}
0 &\equiv d(d\theta_1) \mod{\{\theta_1, \om^1, \om^2\}} && \\
&\equiv -E_1\, \theta_2 \w \om^3 \w \om^4  & \Rightarrow\  & E_1 = 
0. \\
0 &\equiv d(d\theta_2) \mod{\{\theta_2, \om^3, \om^4\}} && \\
&\equiv -E_2\, \theta_1 \w \om^1 \w \om^2  & \Rightarrow\  & E_2 = 
0. \\
0 &\equiv d(d\theta_1) \mod{\{\theta_1, \om^2\}} && \\
&\equiv -(F_1 + A_1 C_2)\, \om^1 \w \om^3 \w \om^4 \qquad & 
\Rightarrow\  & F_1 = -A_1 C_2. \\
0 & \equiv d(d\theta_1) \mod{\{\theta_1, \om^1\}} && \\
&\equiv (-F_2 + A_1 C_1)\, \om^2 \w \om^3 \w \om^4 \qquad & 
\Rightarrow\  & F_2 = A_1 C_1. \\
0 &\equiv d(d\theta_2) \mod{\{\theta_2, \om^4\}} && \\
&\equiv -(F_3 + A_2 C_4)\, \om^1 \w \om^2 \w \om^3 \qquad & 
\Rightarrow\  & F_3 = -A_2 C_4. \\
0 & \equiv d(d\theta_2) \mod{\{\theta_2, \om^3\}} && \\
&\equiv (-F_4 + A_2 C_3)\, \om^1 \w \om^2 \w \om^4 \qquad & 
\Rightarrow\  & F_4 = A_2 C_3.
\end{alignat*}

Next we examine how the functions $A_i, B_i, C_i$ vary if we change 
from one 0-adapted frame to another.  A computation shows that under a 
transformation of the form \eqref{G0freedom}, we have
\begin{gather*}
\tilde{A}_1 = \frac{(a_{11} a_{22} - a_{12} a_{21})}{(b_{11} b_{22} - 
b_{12} b_{21})} A_1 \\[0.1in]
\tilde{A}_2 = \frac{(b_{11} b_{22} - b_{12} b_{21})}{(a_{11} a_{22} - 
a_{12} a_{21})} A_2 
\end{gather*}
\begin{gather*}
\begin{bmatrix} \tilde{B}_1 \\[0.1in] \tilde{B}_2 \end{bmatrix} = 
(b_{11} b_{22} - b_{12} b_{21}) \begin{bmatrix} a_{22} & -a_{12} \\[0.1in] 
-a_{21} & a_{11} \end{bmatrix} \begin{bmatrix} B_1 \\[0.1in] B_2 
\end{bmatrix}\\[0.1in]
\begin{bmatrix} \tilde{B}_3 \\[0.1in] \tilde{B}_4 \end{bmatrix} = 
(a_{11} a_{22} - a_{12} a_{21}) \begin{bmatrix} b_{22} & -b_{12} \\[0.1in] 
-b_{21} & b_{11} \end{bmatrix} \begin{bmatrix} B_3 \\[0.1in] B_4 
\end{bmatrix}
\end{gather*}
\begin{gather*}
\begin{bmatrix} \tilde{C}_1 \\[0.1in] \tilde{C}_2 \end{bmatrix} = 
\frac{(b_{11} b_{22} - b_{12} b_{21})}{(a_{11} a_{22} - a_{12} a_{21})} 
\begin{bmatrix} a_{22} & -a_{12} \\[0.1in] 
-a_{21} & a_{11} \end{bmatrix} \begin{bmatrix} C_1 \\[0.1in] C_2 
\end{bmatrix}\\[0.1in]
\begin{bmatrix} \tilde{C}_3 \\[0.1in] \tilde{C}_4 \end{bmatrix} = 
\frac{(a_{11} a_{22} - a_{12} a_{21})}{(b_{11} b_{22} - b_{12} b_{21})}
\begin{bmatrix} b_{22} & -b_{12} \\[0.1in] 
-b_{21} & b_{11} \end{bmatrix} \begin{bmatrix} C_3 \\[0.1in] C_4 
\end{bmatrix} .
\end{gather*}
From this we see that the functions $A_1, A_2$ and the vectors 
\[ \begin{bmatrix} B_1 \\[0.1in] B_2 \end{bmatrix}, \qquad
\begin{bmatrix} B_3 \\[0.1in] B_4 \end{bmatrix},\qquad
\begin{bmatrix} C_1 \\[0.1in] C_2 \end{bmatrix},\qquad
\begin{bmatrix} C_3 \\[0.1in] C_4 \end{bmatrix} \]
are {\em relative invariants}: if they vanish for any 0-adapted 
coframing, then they vanish for every 0-adapted coframing.  

The general procedure in the method of equivalence is to choose a 
0-adapted coframing that normalizes the torsion terms as much 
as possible.  This has the effect of reducing the structure group to a 
subgroup $G_1 \subset G_0$ which preserves the normalized
torsion terms.  This in turn introduces new torsion terms, which can 
then be further normalized, etc.  Ideally, this process eventually 
leads to a uniquely determined coframing whose torsion terms are 
invariants of the system $\calJ$ on $\calB$.  Even in those cases 
where a unique coframing is not obtained, it may be possible to 
reduce the structure group to the point that some of the torsion 
terms are uniquely determined.  Our hypothesis that $\calB$ is 
homogeneous implies that once the structure group has been reduced to 
the point that it acts trivially on 
a torsion term, that term must be constant on 
$\calB$.

In order 
to proceed with the method of equivalence, we will divide into cases 
depending on whether certain of these invariants are zero or nonzero.

\section{Case 1: $[C_1 \ \ C_2] = [C_3 \ \ C_4] = [0 \ \ 0]$}

Suppose that $C_1 = C_2 = C_3 = C_4 = 0$.  Differentiating the 
structure equations yields
\begin{alignat*}{2}
0 &\equiv d(d\om^1) \mod{\{\theta_1, \om^1, \om^2\}} && \\
&\equiv B_1\, \theta_2 \w \om^3 \w \om^4  & \Rightarrow\  & B_1 = 
0. \\
0 &\equiv d(d\om^2) \mod{\{\theta_1, \om^1, \om^2\}} && \\
&\equiv B_2\, \theta_2 \w \om^3 \w \om^4  & \Rightarrow\  & B_2 = 
0. \\
0 &\equiv d(d\om^3) \mod{\{\theta_2, \om^3, \om^4\}} && \\
&\equiv B_3\, \theta_1 \w \om^1 \w \om^2  & \Rightarrow\  & B_3 = 
0. \\
0 &\equiv d(d\om^4) \mod{\{\theta_2, \om^3, \om^4\}} && \\
&\equiv B_4\, \theta_1 \w \om^1 \w \om^2  & \Rightarrow\  & B_4 = 
0. 
\end{alignat*}
Now we see from the structure equations that
\begin{align*}
& \left. \begin{array}{l}
d\om^1 \equiv 0 \\[0.1in]
d\om^2 \equiv 0 
\end{array} \right\} \mod{\{\om^1, \om^2\}} \\[0.1in]
& \left. \begin{array}{l}
d\om^3 \equiv 0 \\[0.1in]
d\om^4 \equiv 0
\end{array} \right\} \mod{\{\om^3, \om^4\}}.
\end{align*}
Therefore the systems $\{\om^1, \om^2\}$ and $\{\om^3, 
\om^4\}$ are completely integrable; this implies that there exist functions 
$X, Y, P, Q$ on $\calB$ (in fact, these functions are well-defined 
on $M_1$ and $M_2$) such that
\begin{align*}
\{\om^1, \om^2\} & = \{dX, dP\} \\
\{\om^3, \om^4\} & = \{dY, dQ\}.
\end{align*}
The forms $\{\theta_1, dX, 
dY, dP, dQ\}$ comprise a coframing on $M_1$, and by scaling
$\theta_1$ if necessary, we can assume that
\begin{equation}
d\theta_1 = (R_1\, dX + R_2\, dY + R_3\, dP + R_4\, dQ) \w 
\theta_1 + A\, dX \w dP + dY \w dQ \label{case1dtheta}
\end{equation}
for some functions $A, R_i$ on $M_1$ with $A \neq 0$.  Differentiating 
\eqref{case1dtheta}
and reducing modulo $\theta_1$ and $dX \w dP$ yields
\[ (R_1\, dX + R_3\, dP) \w dY \w dQ \equiv 0 \mod{\{\theta_1, dX \w dP\}} 
\qquad \Rightarrow\  R_1 = R_3 = 0. \]
Now differentiating \eqref{case1dtheta} and reducing modulo $\theta_1$ 
yields
\[ [dA - A(R_2\, dY + R_4\, dQ)] \w dX \w dP \equiv 0 \mod{\theta_1}. 
\]
Therefore there exist functions $A_0, A_1, A_3$ on $M_1$ such that
\[ dA = A_0\, \theta_1 + A_1\, dX + A R_2\, dY + A_3\, dP + A R_4\, 
dQ. \]
Now differentiating \eqref{case1dtheta} yields
\[ (dR_2 \w dY + dR_4 \w dQ + A_0\, dX \w dP) \w \theta_1 = 0. \]
Reducing modulo $\{dY, dQ\}$ shows that $A_0 = 0$, and so this equation 
becomes
\[ (dR_2 \w dY + dR_4 \w dQ) \w \theta_1 = 0. \]
By Cartan's lemma (see \cite{BCG3} for details), this implies that
there exist functions $R_{20}, R_{40}, 
R_{22}$, $R_{24}, R_{44}$ such that
\begin{align}
dR_2 & = R_{20}\, \theta_1 + R_{22}\, dY + R_{24}\, dQ  \label{case1Req}\\
dR_4 & = R_{40}\, \theta_1 + R_{24}\, dY + R_{44}\, dQ . \notag
\end{align}
Differentiating these equations and reducing modulo $\{\theta_1, dY, 
dQ\}$ yields
\[ \left. \begin{array}{l}
0 = AR_{20}\, dX \w dP \\
0 = AR_{40}\, dX \w dP 
\end{array} \right\} \mod{\{\theta_1, dY, dQ\}}.  \]
Therefore $R_{20} = R_{40} = 0$, and from \eqref{case1Req} we see that
\[ d(R_2\, dY + R_4\, dQ) = 0. \]

Let $\lam(Y, Q)$ be a nonvanishing function such that
\[ \lam^{-1}\, d\lam = R_2\, dY + R_4\, dQ \]
and let $\tilde{\theta}_1 = \lam^{-1} \theta_1.$  Then
\begin{equation} 
d\tilde{\theta}_1 = \lam^{-1}A\, dX \w dP + \lam^{-1}\, dY \w dQ. 
\label{dthetatilde}
\end{equation}
Differentiating yields 
\[ d(\lam^{-1}A) \w dX \w dP = 0 \]
and so the function $\tilde{A} = \lam^{-1}A$ is a function of $X$ and 
$P$ alone.  It follows that 
\[ d(\lam^{-1}A\, dX \w dP) = d(\lam^{-1}\, dY \w dQ) = 0. \]
Therefore, by Darboux's Theorem there 
exist new functions
\begin{align*}
x_1 & = x_1(X, P) \\
p_1 & = p_1(X, P) \\
y_1 & = y_1(Y, Q) \\
q_1 & = q_1(Y, Q)
\end{align*}
such that equation \eqref{dthetatilde} takes the form
\[ d\tilde{\theta}_1 = -dp_1 \w dx_1 - dq_1 \w dy_1 \] 
and
\[ \{\om^1, \om^2 \} = \{dx_1, dp_1\} , \qquad
\{\om^3, \om^4 \} = \{dy_1, dq_1\} . \]
Now by Pfaff's Theorem there must exist a function $z_1$ on $M_1$ such that
\[ \tilde{\theta}_1 = dz_1 - p_1\, dx_1 - q_1\, dy_1 \]
and we see that $\calI_1$ is the ideal corresponding to the wave 
equation
\[ z_{xy} = 0. \]

By the same argument, the ideal $\calI_2$ on $M_2$ also represents the 
wave equation,
and the B\"acklund transformation is given by equations of the form
\begin{align*}
x_2 & = x_2(x_1, p_1) \\
p_2 & = p_2(x_1, p_1) \\
y_2 & = y_2(y_1, q_1) \\
q_2 & = q_2(y_1, q_1).
\end{align*}
These may be written in PDE notation as
\begin{align*}
\bar{x} & = \bar{x}(x, z_x) \\
\bar{z}_{\bar{x}} & = \bar{z}_{\bar{x}}(x, z_x) \\
\bar{y} & = \bar{y}(y, z_y) \\
\bar{z}_{\bar{y}} & = \bar{z}_{\bar{y}}(y, z_y),
\end{align*}
and the nondegeneracy conditions imply that
\[ 0 \neq \frac{\partial p_2}{\partial p_1} \frac{\partial x_2}{\partial 
x_1} - \frac{\partial p_2}{\partial x_1} \frac{\partial x_2}{\partial 
p_1} \neq \frac{\partial q_2}{\partial q_1} \frac{\partial y_2}{\partial 
y_1} - \frac{\partial q_2}{\partial y_1} \frac{\partial y_2}{\partial 
q_1} \neq 0. \]
These transformations are somewhat more general than typical
point transformations (or even gauge transformations), in that they do 
not necessarily preserve the space of independent variables.

Thus we have the following theorem.

\begin{theorem}\label{case1thm}
Let $\calB \subset M_1 \times M_2$ be a B\"acklund transformation with 
$C_1 = C_2 = C_3 = C_4 = 0$.  Then $\calB$ is locally contact 
equivalent to a transformation between solutions of the wave equation
\[ z_{xy} = 0 \]
with the property that given any solution, the new 
solutions given by the transformation may be obtained by quadrature.
\end{theorem}
Note that in this case the assumption of homogeneity was not 
necessary.  In the remaining cases, however, homogeneity will play a 
crucial role in the analysis.

\section{Case 2: $[C_1 \ \ C_2] = [0 \ \ 0], \ [C_3 \ \ C_4] \neq [0 \ \ 
0]$}

Suppose that $C_1 = C_2 = 0$, but $C_3$ and $C_4$ are not both zero.  
By a transformation of the form \eqref{G0freedom}, we can arrange 
that 
\[ C_3 = 0,\ \  C_4 = 1,\ \  A_1 = 1. \]
A coframing satisfying this condition will be called {\em 1-adapted}.  
If $\{\theta_1,\, \theta_2,\, \om^1,\, 
\om^2,\, \om^3$, $\om^4\}$ is a 1-adapted coframing, then any other 
1-adapted coframing $\{\tilde{\theta}_1,\, \tilde{\theta}_2,\, \tilde{\om}^1,\, 
\tilde{\om}^2$, 
$\tilde{\om}^3,\, \tilde{\om}^4\}$ has the form
\begin{equation}
 \begin{bmatrix} \tilde{\theta}_1\\[0.1in] \tilde{\theta}_2\\[0.1in] 
\tilde{\om}^1\\[0.1in] \tilde{\om}^2\\[0.1in] \tilde{\om}^3\\[0.1in] 
\tilde{\om}^4 
\end{bmatrix} =  
\begin{bmatrix} \scriptstyle{a_{11}a_{22} - a_{12}a_{21}} & 0 & 
0 & 0 & 0 & 0 \\[0.1in] 0 & \scriptstyle{a_{11}a_{22} - a_{12}a_{21}} & 
0 & 0 & 0 & 0 
\\[0.1in] 0 & 0 & a_{11} & a_{12} & 0 & 0 
\\[0.1in] 0 & 0 & a_{21} & a_{22} & 0 & 0 
\\[0.1in] 0 & 0 & 0 & 0 & 1 & 0 \\[0.1in] 0 & 0 & 0 & 
0 & b_{21} & \scriptstyle{a_{11}a_{22} - a_{12}a_{21}}
\end{bmatrix}^{-1} \begin{bmatrix} \theta_1\\[0.1in] \theta_2\\[0.1in] 
\om^1\\[0.1in] 
\om^2\\[0.1in] \om^3\\[0.1in] \om^4 \end{bmatrix}.  \label{case2G1freedom}
\end{equation}

The same computation as in the previous section shows that $B_1 = 
B_2 = 0$.  Furthermore, we have
\begin{align*}
0 & \equiv d(d\om^3) \mod{\{\theta_1, \om^3, \om^4\}} \\
& \equiv (\bet_2 + B_3\, \theta_2) \w \om^1 \w \om^2. \\
0 & \equiv d(d\om^3) \mod{\{\theta_2, \om^3, \om^4\}} \\
& \equiv (\bet_2 - B_3\, \theta_1) \w \om^1 \w \om^2. 
\end{align*}
Together, these equations imply that
\[ \bet_2 = B_3\, \theta_1 - B_3\, \theta_2 + H_1\, \om^1 + H_2\, 
\om^2 + H_3\, \om^3 + H_4\, \om^4 \]
for some functions $H_1, H_2, H_3, H_4$ on $\calB$.  Similarly, computing
$d(d\om^4) \equiv 0$ modulo $\{\theta_1, \om^3, \om^4\}$ and 
$\{\theta_2, \om^3, \om^4\}$ shows that
\[ \bet_4 = \alp_1 + \alp_4 + B_4\, \theta_1 - B_4\, \theta_2 + J_1\, 
\om^1 + J_2\, \om^2 + J_3\, \om^3 + J_4\, \om^4 \]
for some functions $J_1, J_2, J_3, J_4$ on $\calB$.  By taking 
advantage of the ambiguity 
\eqref{G0connfreedom} in the forms $\bet_i$, we can assume that
\[ H_3 = H_4 = J_3 = 0. \]
Now computing $d(d\theta_1) \equiv 0 \mod{\theta_1}$ shows that
\[ \bet_1 = K_0\, \theta_1 + B_4\, \theta_2 + K_1\, \om^1 + K_2\, 
\om^2 + \om^3 - J_4\, \om^4 \]
for some functions $K_0, K_1, K_2$ on $\calB$.

Under a transformation of the form \eqref{case2G1freedom}, the 
function $A_2$ remains unchanged; therefore by our assumption of 
homogeneity it must be constant.  Moreover, the nondegeneracy 
assumptions imply that $A_2 \neq 0, 1$.  So 
\begin{align*}
0 & \equiv d(d\theta_2) \mod{\theta_2} \\
& \equiv -A_2 [(K_0 + B_4)\, \theta_1 + (K_1 + J_1)\, \om^1 + (K_2 + 
J_2)\, \om^2] \w \om^3 \w \om^4,
\end{align*}
which implies that
\begin{align*}
K_0 & = -B_4 \\
K_1 & = -J_1 \\
K_2 & = -J_2.
\end{align*}
Now we have 
\begin{align*}
0 & = d(d\theta_1) = \Upsilon_1 \w \theta_1  \\
0 & = d(d\theta_2) = \Upsilon_2 \w \theta_2 
\end{align*}
where 
\begin{align*}
\Upsilon_1 & = d\alp_1 + d\alp_4 + (J_1\, \om^1 + J_2\, \om^2) \w 
\om^3 - (H_1\, \om^1 + H_2\, \om^2) \w \om^4 - J_4\, \om^3 \w \om^4 \\
\Upsilon_2 & = d\alp_1 + d\alp_4 + A_2 [(J_1\, \om^1 + J_2\, \om^2) \w 
\om^3 - (H_1\, \om^1 + H_2\, \om^2) \w \om^4 - J_4\, \om^3 \w \om^4].
\end{align*}
These equations imply that $\Upsilon_1$ must be a multiple of $\theta_1$ and 
$\Upsilon_2$ must be a multiple of $\theta_2$, so 
\begin{align*}
0 & \equiv \Upsilon_2 - \Upsilon_1 \mod{\{\theta_1, \theta_2\}} \\
& \equiv (A_2 - 1)[(J_1\, \om^1 + J_2\, \om^2) \w 
\om^3 - (H_1\, \om^1 + H_2\, \om^2) \w \om^4 - J_4\, \om^3 \w \om^4].
\end{align*}
Therefore,
\[ H_1 = H_2 = J_1 = J_2 = J_4 = 0. \]

The structure equations for a 1-adapted coframing now take the form
\begin{equation}
\begin{bmatrix} d\theta_1\\[0.1in] d\theta_2\\[0.1in] d\om^1\\[0.1in] 
d\om^2\\[0.1in] d\om^3\\[0.1in] d\om^4
\end{bmatrix} = 
-\begin{bmatrix}\alp_1 + \alp_4 & 0 & 0 & 0 & 0 & 0 \\[0.1in] 0 & \alp_1 + 
\alp_4 & 0 & 0 & 0 & 0 \\[0.1in] 0 & 0 & \alp_1 & \alp_2 & 0 & 0 
\\[0.1in] 0 & 0 & 
\alp_3 & \alp_4 & 0 & 0 \\[0.1in] 0 & 0 & 0 & 0 & 0 & 0
\\[0.1in] 0 & 0 & 0 & 
0 & \bet_3 & \alp_1 + \alp_4 \end{bmatrix} \w
\begin{bmatrix} \theta_1\\[0.1in] \theta_2\\[0.1in] \om^1\\[0.1in] 
\om^2\\[0.1in] \om^3\\[0.1in] \om^4 
\end{bmatrix} + 
\begin{bmatrix} \Theta_1\\[0.1in] \Theta_2\\[0.1in] \Omega^1\\[0.1in] 
\Omega^2\\[0.1in] \Omega^3\\[0.1in] \Omega^4 
\end{bmatrix} \label{case2G1struct}
\end{equation}
where 
\begin{align*}
\Theta_1 & = \om^1 \w \om^2 + \om^3 \w (\om^4 - \theta_1)  
\\
\Theta_2 & = \om^1 \w \om^2 + A_2\, \om^3 \w (\om^4 - \theta_2) \\
\Omega^1 & = 0  \\
\Omega^2 & = 0  \\
\Omega^3 & = B_3\, (\theta_1 - \om^4) \w (\theta_2 - \om^4) + B_4\, 
(\theta_1 - \theta_2) \w \om^3 \\
\Omega^4 & = B_4\,(\theta_1 - \om^4) \w (\theta_2 - \om^4) + \om^1 \w 
\om^2.
\end{align*}

Now 
\[ 0 = d(d\theta_1) = -d(\alp_1 + \alp_4) \w \theta_1, \]
and so 
\[ d(\alp_1 + \alp_4) = \psi \w \theta_1 \]
for some 1-form $\psi$.  Differentiating this equation and reducing 
modulo $\theta_1$ yields
\[ -\psi \w (\om^1 \w \om^2 + \om^3 \w \om^4) \equiv 0 
\mod{\theta_1}. \]
But since $\om^1 \w \om^2 + \om^3 \w \om^4$ is nondecomposable, this 
implies that
\[ \psi \equiv 0 \mod{\theta_1} \]
and hence that
\[ d(\alp_1 + \alp_4) = 0. \]
Therefore, there exists a nonvanishing function $\lambda$ on $\calB$ 
such that
\[ \lam^{-1}\, d\lam = -(\alp_1 + \alp_4). \]

We can choose a new 1-adapted coframing in which $\theta_1$ is 
replaced by $\lam^{-1} \theta_1$.
This coframing will have the property that 
\[ \alp_1 + \alp_4 = 0. \]
Now we have
\[ d(\om^1 \w \om^2) = 0, \]
and so by Darboux's Theorem there exist functions $x, p$ on $\calB$ (which 
are also well-defined on $M_1$ and $M_2$) such that
\[ \om^1 \w \om^2 = dx \w dp. \]
Therefore
\begin{align*}
0 & = d(d\theta_1) = d(\om^3 \w (\om^4 - \theta_1)) \\
0 & = d(d\theta_2) = d(A_2\, \om^3 \w (\om^4 - \theta_2)).
\end{align*}
Again by Darboux's Theorem, there exist functions $y_1, q_1$ on $M_1$ and 
$y_2, q_2$ on $M_2$ such that
\begin{align*}
\om^3 \w (\om^4 - \theta_1) & = dy_1 \w dq_1 \\
A_2\, \om^3 \w (\om^4 - \theta_2) & = dy_2 \w dq_2.
\end{align*}
By Pfaff's Theorem there exist functions $z_1$ on $M_1$ and $z_2$ on 
$M_2$ such that
\begin{align*}
\theta_1 & = dz_1 - p\, dx - q_1\, dy_1 \\
\theta_2 & = dz_2 - p\, dx - q_2\, dy_2.
\end{align*}
The ideals $\calI_1, \calI_2$ now take the form
\begin{gather*}
\calI_1 = \{\theta_1, \, \om^1 \w \om^2,\, \om^3 \w (\om^4 - 
\theta_1) \} = \{ dz_1 - p\, dx - q_1\, dy_1,\, dx \w dp,\, dy_1 \w 
dq_1\} \\
\calI_2 = \{\theta_2, \, \om^1 \w \om^2,\, A_2\, \om^3 \w (\om^4 - 
\theta_2) \} = \{ dz_2 - p\, dx - q_2\, dy_2,\, dx \w dp,\, dy_2 \w 
dq_2\}.
\end{gather*}
Both represent the wave equation 
\[ z_{xy} = 0, \]
and the B\"acklund transformation is given by equations of the form
\begin{align*}
x_2 & = x_1 \\
y_2 & = y_2(x_1, y_1, z_1, z_2, p_1, q_1) \\
p_2 & = p_1 \\
q_2 & = q_2(x_1, y_1, z_1, z_2, p_1, q_1)),
\end{align*}
or, in PDE notation,
\begin{align*}
\bar{x} & = x \\
\bar{y} & = \bar{y}(x, y, z, \bar{z}, z_x, z_y) \\
\bar{z}_{\bar{x}} & = z_x \\
\bar{z}_{\bar{y}} & = \bar{z}_{\bar{y}}(x, y, z, \bar{z}, z_x, z_y).
\end{align*}
As in the previous case, these transformations do not in general preserve 
the space of independent variables.

Thus we have the following theorem.

\begin{theorem}\label{case2thm}
Let $\calB \subset M_1 \times M_2$ be a homogeneous B\"acklund 
transformation with one of the vectors $[C_1 \ \ C_2], \ [C_3 \ \ 
C_4]$ identically zero and the other nonzero.
Then $\calB$ is locally contact 
equivalent to a transformation between solutions of the wave equation
\[ z_{xy} = 0. \]
\end{theorem}

\section{Case 3: $[C_1 \ \ C_2], \ [C_3 \ \ C_4] \neq [0 \ \ 
0]$}

Suppose that the vectors $[C_1 \ \ C_2], \ [C_3 \ \ C_4]$ are both 
nonzero.  By a transformation of the form \eqref{G0freedom}, we can arrange 
that 
\[ C_1 = C_3 = 0, \qquad C_2 = C_4 = 1. \]
A coframing satisfying this 
condition is called {\em 1-adapted}.  If $\{\theta_1,\, \theta_2,\, \om^1,\, 
\om^2,\, \om^3$, $\om^4\}$ is a 1-adapted coframing, then any other 
1-adapted coframing $\{\tilde{\theta}_1,\, \tilde{\theta}_2,\, \tilde{\om}^1,\, 
\tilde{\om}^2,\, 
\tilde{\om}^3$, $\tilde{\om}^4\}$ has the form
\begin{equation}
 \begin{bmatrix} \tilde{\theta}_1\\[0.1in] \tilde{\theta}_2\\[0.1in] 
\tilde{\om}^1\\[0.1in] \tilde{\om}^2\\[0.1in] \tilde{\om}^3\\[0.1in] 
\tilde{\om}^4 
\end{bmatrix} =  
\begin{bmatrix} a_{22} & 0 & 
0 & 0 & 0 & 0 \\[0.1in] 0 &  b_{22} & 0 & 0 & 0 & 0 
\\[0.1in] 0 & 0 & \frac{b_{22}}{a_{22}} & 0 & 0 & 0 \\[0.1in] 0 & 0 & a_{21} & 
a_{22} & 0 & 0 
\\[0.1in] 0 & 0 & 0 & 0 & \frac{a_{22}}{b_{22}} & 0 \\[0.1in] 0 & 0 & 0 & 
0 & b_{21} & b_{22}
\end{bmatrix}^{-1} \begin{bmatrix} \theta_1\\[0.1in] \theta_2\\[0.1in] 
\om^1\\[0.1in] 
\om^2\\[0.1in] \om^3\\[0.1in] \om^4 \end{bmatrix}.  \label{case3G1freedom}
\end{equation}
Similarly to the previous case, computing 
\begin{align*}
d(d\om^1) & \equiv 0 \mod{\{\theta_1, \om^1, \om^2\}} \\
d(d\om^1) & \equiv 0 \mod{\{\theta_2, \om^1, \om^2\}} \\
d(d\om^3) & \equiv 0 \mod{\{\theta_1, \om^3, \om^4\}} \\
d(d\om^3) & \equiv 0 \mod{\{\theta_2, \om^3, \om^4\}}
\end{align*}
shows that 
\begin{align*}
\alp_2 & = A_2 B_1\, \theta_1 - B_1\, 
\theta_2 + G_1\, \om^1 + G_2\, \om^2 + G_3\, \om^3 + G_4\, \om^4 \\
\bet_2 & = B_3\, \theta_1 - A_1 B_3\, 
\theta_2 + H_1\, \om^1 + H_2\, \om^2 + H_3\, \om^3 + H_4\, \om^4
\end{align*}
for some functions $G_i, H_i$ on $\calB$.  By taking advantage of the 
ambiguity \eqref{G0connfreedom}, we can assume that
\[ G_1 = G_2 = H_3 = H_4 = 0. \]
Now computing
\begin{align*}
d(d\om^2) & \equiv 0 \mod{\{\theta_1, \om^1, \om^2\}} \\
d(d\om^2) & \equiv 0 \mod{\{\theta_2, \om^1, \om^2\}} \\
d(d\om^4) & \equiv 0 \mod{\{\theta_1, \om^3, \om^4\}} \\
d(d\om^4) & \equiv 0 \mod{\{\theta_2, \om^3, \om^4\}}
\end{align*}
shows that 
\begin{align*}
\alp_1 & = \bet_4 - \alp_4 - B_4\, \theta_1 + A_1 B_4\, 
\theta_2 + J_1\, \om^1 + J_2\, \om^2 + J_3\, \om^3 + J_4\, \om^4 \\
\bet_1 & = \alp_4 - \bet_4 - A_2 B_2\, \theta_1 + B_2\, 
\theta_2 + K_1\, \om^1 + K_2\, \om^2 + K_3\, \om^3 + K_4\, \om^4
\end{align*}
for some functions $J_i, K_i$ on $\calB$.  Using some of the remaining 
ambiguity \eqref{G0connfreedom}, we can assume that 
\[ J_3 = K_1 = 0. \]
The structure equations for a 1-adapted coframing now take the form
\begin{equation}
\begin{bmatrix} d\theta_1\\[0.1in] d\theta_2\\[0.1in] d\om^1\\[0.1in] 
d\om^2\\[0.1in] d\om^3\\[0.1in] d\om^4
\end{bmatrix} = 
-\begin{bmatrix}\alp_4 & 0 & 0 & 0 & 0 & 0 \\[0.1in] 0 & 
\bet_4 & 0 & 0 & 0 & 0 \\[0.1in] 0 & 0 & \bet_4 - \alp_4 & 0 & 0 & 0 
\\[0.1in] 0 & 0 & 
\alp_3 & \alp_4 & 0 & 0 \\[0.1in] 0 & 0 & 0 & 0 & \alp_4 - \bet_4 & 0
\\[0.1in] 0 & 0 & 0 & 
0 & \bet_3 & \bet_4 \end{bmatrix} \w
\begin{bmatrix} \theta_1\\[0.1in] \theta_2\\[0.1in] \om^1\\[0.1in] 
\om^2\\[0.1in] \om^3\\[0.1in] \om^4 
\end{bmatrix} + 
\begin{bmatrix} \Theta_1\\[0.1in] \Theta_2\\[0.1in] \Omega^1\\[0.1in] 
\Omega^2\\[0.1in] \Omega^3\\[0.1in] \Omega^4 
\end{bmatrix} \label{case3G1struct}
\end{equation}
where 
\begin{align*}
\Theta_1 & = \theta_1 \w (B_2\, \theta_2 + A_1\, \om^1 + K_2\, 
\om^2 + K_3\, \om^3 + K_4\, \om^4) + A_1\, \om^1 \w \om^2 + \om^3 \w \om^4
\\
\Theta_2 & = \theta_2 \w (-B_4\, \theta_1 + J_1\, \om^1 + J_2\, 
\om^2 + A_2\, \om^3 + J_4\, \om^4) + \om^1 \w \om^2 + A_2\, \om^3 \w \om^4  \\
\Omega^1 & = \om^1 \w (-B_4\, \theta_1 + A_1 B_4\, \theta_2 + J_2\, 
\om^2 + J_4\, \om^4)\\
& \qquad \qquad  + \om^2 \w (A_2 B_1\, \theta_1 - B_1\, \theta_2 + 
G_3\, \om^3 + G_4\, \om^4) + B_1\, \theta_1 \w \theta_2  \\
\Omega^2 & = B_2\, \theta_1 \w \theta_2 + \om^3 \w \om^4  \\
\Omega^3 & = \om^3 \w (-A_2 B_2\, \theta_1 + B_2\, \theta_2 + K_2\, 
\om^2 + K_4\, \om^4)  \\
& \qquad \qquad + \om^4 \w (B_3\, \theta_1 - A_1 B_3\, \theta_2 + 
H_1\, \om^1 + H_2\, \om^2) + B_3\, \theta_1 \w \theta_2 \\
\Omega^4 & = B_4\, \theta_1 \w \theta_2 + \om^1 \w \om^2  .
\end{align*}

A computation shows that under a transformation of the form 
\eqref{case3G1freedom}, we have
\begin{gather*}
\begin{bmatrix} \tilde{B}_1 \\[0.1in] \tilde{B}_2 \end{bmatrix} = 
\begin{bmatrix} (a_{22})^2 & 0 \\[0.1in] 
-a_{21} a_{22} & b_{22} \end{bmatrix} \begin{bmatrix} B_1 \\[0.1in] B_2 
\end{bmatrix}\\[0.1in]
\begin{bmatrix} \tilde{B}_3 \\[0.1in] \tilde{B}_4 \end{bmatrix} = 
\begin{bmatrix} (b_{22})^2 & 0 \\[0.1in] 
-b_{21} b_{22} & a_{22} \end{bmatrix} \begin{bmatrix} B_3 \\[0.1in] B_4 
\end{bmatrix}.
\end{gather*}
In particular, the functions $B_1, B_3$ are now relative invariants.
In order to proceed further, we will need to divide into cases depending 
on the values of the $B_i$.  First we prove the following lemma:

\begin{lemma}\label{Bveclemma}
For any 1-adapted coframing as above, the vectors $[B_1 \ \ B_2], \ 
[B_3 \ \ B_4]$ are either both zero or both nonzero.
\end{lemma}

\begin{pf}
Suppose that $B_1 = B_2 = 0.$  Then 
\begin{align*}
0 & \equiv d(d\om^2) \mod{\{\om^1, \om^2\}} \\
& \equiv \theta_1 \w \theta_2 \w (B_3\, \om^4 - B_4\, \om^3).
\end{align*}
Therefore, $B_3 = B_4 = 0.$  A similar argument demonstrates the 
converse.
\end{pf}

\section{Case 3A: $[B_1 \ \ B_2] = [B_3 \ \ B_4] = [0 \ \ 0]$}

\begin{proposition}
In this case we can choose a 1-adapted coframing for which $A_1, A_2$ 
are constant and $\{\alp_4, \bet_4\} \in {\rm 
span}\{\om^1, \om^2, \om^3, \om^4\}$.
\end{proposition}

\begin{pf}
Under a transformation of the form \eqref{case3G1freedom}, we have
\begin{align*}
\tilde{A}_1 & = \frac{b_{22}}{a_{22}} A_1 \\
\tilde{A}_2 & = \frac{a_{22}}{b_{22}} A_2 \\
\tilde{G}_3 & = \frac{a_{22}^3}{b_{22}^2} G_3 + \frac{a_{22}^2 
b_{21}}{b_{22}} G_4 \\
\tilde{G}_4 & = a_{22}^2\, G_4 \\
\tilde{H}_1 & = \frac{b_{22}^3}{a_{22}^2} H_1 + \frac{b_{22}^2 
a_{21}}{a_{22}} H_2 \\
\tilde{H}_2 & = b_{22}^2\, H_2\\
\tilde{J}_1 & = \frac{b_{22}}{a_{22}} J_1 + a_{21}\, J_2 + 
\frac{b_{22}^2 b_{21}}{a_{22}^2} H_1 + \frac{b_{22} a_{21} 
b_{21}}{a_{22}} H_2 \\
\tilde{J}_2 & = a_{22}\, J_2 + b_{22} b_{21}\, H_2 \\
\tilde{J}_4 & = b_{22}\, J_4 \\
\tilde{K}_2 & = a_{22}\, K_2 \\
\tilde{K}_3 & = \frac{a_{22}}{b_{22}} K_3 + b_{21}\, K_4 + 
\frac{a_{22}^2 a_{21}}{b_{22}^2} G_3 + \frac{a_{22} a_{21} 
b_{21}}{b_{22}} G_4 \\
\tilde{K}_4 & = b_{22}\, K_4 + a_{22} a_{21}\, G_4
\end{align*}
By the homogeneity assumption, we can choose a 1-adapted coframing for 
which all the torsion functions are constants.  For such a coframing, 
we have
\begin{align*}
0 & \equiv d(d\theta_1) \mod{\theta_1}  \\
& \equiv [A_1\, (\alp_4 - \bet_4) + (A_1 K_3 - 1)\, \om^3 + (A_1 
K_4 - A_1 J_4)\, \om^4] \w \om^1 \w \om^2 \\
0 & \equiv d(d\theta_2) \mod{\theta_2}  \\
& \equiv [A_2\, (\bet_4 - \alp_4) + (A_2 J_1 - 1)\, \om^1 + (A_2 
J_2 - A_2 K_2)\, \om^2] \w \om^3 \w \om^4.
\end{align*}
Together these equations imply that
\[ \bet_4 = \alp_4 + \left( \frac{1}{A_2} - J_1 \right)\, \om^1 + 
(K_2 - J_2)\, \om^2 + \left( K_3 - \frac{1}{A_1} \right)\, \om^3 + 
(K_4 - J_4)\, \om^4, \]
so it suffices to show that $\alp_4 \in {\rm 
span}\{\om^1, \om^2, \om^3, \om^4\}$.  For this we divide into cases 
depending on whether certain of the torsion functions are zero or 
nonzero.
\begin{itemize}
\item{Observe that the functions $G_4, H_2, J_4, K_2$ are relative 
invariants.  Moreover, we have
\begin{align*}
0 & \equiv d(d\om^1) \mod{\{\om^1, \om^3\}} \\
& \equiv -2 G_4\, \alp_4 \w \om^2 \w \om^4 \\
0 & \equiv d(d\om^3) \mod{\{\om^1, \om^3\}} \\
& \equiv 2 H_2\, \alp_4 \w \om^2 \w \om^4 \\
0 & \equiv d(d\om^3) \mod{\{\om^1, \om^2\}} \\
& \equiv -J_4\, \alp_4 \w \om^3 \w \om^4 \\
0 & \equiv d(d\om^1) \mod{\{\om^3, \om^4\}} \\
& \equiv -K_2\, \alp_4 \w \om^1 \w \om^2.
\end{align*}
So if any of these invariants is nonzero, we conclude that $\alp_4 \in {\rm 
span}\{\om^1, \om^2,$ $\om^3, \om^4\}$, as desired.}
\item{If $G_4 = H_2 = J_4 = K_2 = 0$, then the functions $G_3, H_1, 
J_2, K_4$ become relative invariants and we have
\begin{align*}
0 & \equiv d(d\om^1) \mod{\om^1} \\
& \equiv (-G_3\, \alp_4 + G_3 K_4\, \om^4) \w \om^2 \w \om^3 \\
0 & \equiv d(d\om^3) \mod{\om^3} \\
& \equiv (H_1\, \alp_4 - 2 H_1 J_2\, \om^2) \w \om^1 \w \om^4\\
0 & \equiv d(d\om^3) \mod{\{\om^1, \om^4\}} \\
& \equiv J_2\, \alp_4 \w \om^2 \w \om^3 \\
0 & \equiv d(d\om^1) \mod{\{\om^2, \om^3\}} \\
& \equiv -K_4\, \alp_4 \w \om^1 \w \om^4. 
\end{align*}
So if any of these invariants is nonzero, we conclude that $\alp_4 \in {\rm 
span}\{\om^1, \om^2,$ $\om^3, \om^4\}$, as desired.}
\item{If $G_3 = G_4 = H_1 = H_2 = J_2 = J_4 = K_2 = K_4 = 0$, then we 
have 
\begin{align*}
0 & = d(d\theta_1) \\
& = (-d\alp_4 + Q_1\, \om^1 \w \om^3) \w \theta_1 \\
0 & =d(d\theta_2) \\
& = (-d\alp_4 + Q_2\, \om^1 \w \om^3) \w \theta_2
\end{align*}
where $Q_1, Q_2$ are certain functions of the constants $A_1, A_2, 
J_1, K_3$.  Together these equations imply that $Q_1 = Q_2$ and that
\[ d\alp_4 = Q\, \om^1 \w \om^3 + Z\, \theta_1 \w \theta_2 \]
where $Q = Q_1 = Q_2$ is constant and $Z$ is some function on 
$\calB$.  Differentiating this equation and reducing modulo 
$\theta_1$ yields
\begin{align*}
0 & \equiv d(d\alp_4) \mod{\theta_1} \\
& \equiv Z(A_1\, \om^1 \w \om^2 + \om^3 \w \om^4).
\end{align*}
Therefore $Z=0$, and we have
\[ d\alp_4 = Q\, \om^1 \w \om^3. \]
Moreover, since
\begin{gather*}
d\om^1 = \left( K_3 - \frac{1}{A_1} \right) \, \om^1 \w \om^3 \\
d\om^3 = \left( \frac{1}{A_2} - J_1 \right) \, \om^1 \w \om^3
\end{gather*}
we conclude that $d\alp_4$ is a constant multiple of at least one of 
$d\om^1, d\om^3$.  (It is straightforward to check that if $d\om^1 = 
d\om^3 = 0$, then $d\alp_4 = 0$ as well.)  Without loss of generality, 
assume that $d\alp_4$ 
is a constant multiple of $d\om^1$.  Then there exists a nonvanishing 
function $f$ on $\calB$ and a constant $C$ such that
\[ \alp_4 = C\, \om^1 - f^{-1}\,df. \]
The only remaining torsion terms in this case are $A_1, A_2, J_1, K_3$.  
Note that 
these terms all remain unchanged under a transformation of the form 
\eqref{case3G1freedom} with $b_{22} = a_{22}$.  Moreover, if we take 
$a_{22} = b_{22} = f$, the new coframing satisfies the condition 
that
\[ \alp_4 = C\, \om^1 \]
and we have $\alp_4 \in {\rm 
span}\{\om^1, \om^2, \om^3, \om^4\}$, as desired.}
\end{itemize}
\end{pf}

Now suppose that we have chosen a coframing as in the proposition.  
The forms $\{\om^1, \om^2, \om^3, \om^4\}$ form a Frobenius system, 
and so locally there exists a 4-manifold $N$ which is a quotient of 
$\calB$ and for which the 1-forms $\om^1, \om^2, \om^3, \om^4$ are 
semi-basic for the projection $\calB \to N$.  (Here ``locally'' refers 
to the fact that any point in $\calB$ has a neighborhood which 
possesses such 
a quotient, and ``semi-basic'' means that the 
restrictions of the $\om^i$ to the fibers of the projection vanish 
identically.  See \cite{BCG3} for details.)  In fact, this quotient 
fibers through each of the quotients $\pi_i: \calB \to M_i$, as shown.
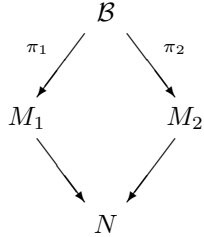
\setlength{\unitlength}{2pt}
\begin{center}
\begin{picture}(40,50)(0,0)
\put(5,25){\makebox(0,0){$M_1$}}
\put(35,25){\makebox(0,0){$M_2$}}
\put(20,45){\makebox(0,0){$\calB$}}
\put(16,41){\vector(-3,-4){9}}
\put(24,41){\vector(3,-4){9}}
\put(7,38){\makebox(0,0){$\scriptstyle{\pi_1}$}}
\put(33,38){\makebox(0,0){$\scriptstyle{\pi_2}$}}
\put(20,5){\makebox(0,0){$N$}}
\put(7,21){\vector(3,-4){9}}
\put(33,21){\vector(-3,-4){9}}
\end{picture}
\end{center}  
Moreover, we have
\begin{align*}
d(\om^1 \w \om^2) & = -\bet_4 \w \om^1 \w \om^2 - J_4\, \om^1 \w 
\om^2 \w \om^4 - \om^1 \w \om^3 \w \om^4 \\
d(\om^3 \w \om^3) & = -\alp_4 \w \om^3 \w \om^4 - K_2\, \om^2 \w 
\om^3 \w \om^4 - \om^1 \w \om^2 \w \om^3,
\end{align*}
and so $d(\om^1 \w \om^2),\ d(\om^3 \w \om^4) \in 
\Lambda^3(\{\om^i\})$.  Therefore the forms $\om^1 \w \om^2,\ \om^3 \w 
\om^4$ are well-defined on $N$.

From the structure equations \eqref{case3G1struct}, we now have 
\[ d\theta_1 = \gamma \w \theta_1 + A_1\, \om^1 \w \om^2 + \om^3 \w 
\om^4 \]
for some 1-form $\gamma$ on $\calB$ which we can assume is a linear 
combination of $\theta_2$ and the $\om^i$; moreover, $A_1$ is constant.  
Differentiating this equation and reducing modulo $\theta_1$ and
$\Lambda^3(\{\om^i\})$ shows that in fact $\gamma$ is a linear 
combination of the $\om^i$ alone.  Then differentiating and reducing 
modulo $\Lambda^3(\{\om^i\})$ yields
\[ 0 \equiv d\gamma \w \theta_1  \mod{\Lambda^3(\{\om^i\})}, \]
but this implies that in fact
\[ d\gamma \w \theta_1 = 0. \]
By an argument identical to that given in Case 2 for the form $\alp_1 + 
\alp_4$, it follows that
\[ d\gamma = 0. \]
Therefore, there exists a nonvanishing function $\lambda$ on $\calB$ 
such that
\[ \gamma = -\lambda^{-1}\, d\lambda. \]
Let $\tilde{\theta}_1 = \lam\, \theta_1$.  Then
\begin{align*}
d\tilde{\theta}_1 & = d\lam \w \theta_1 + \lam\, d\theta_1 \\
& = \lam(A_1\, \om^1 \w \om^2 + \om^3 \w \om^4).
\end{align*}
Since this is a closed form which is semi-basic for the projection 
$\calB \to N$, it is in fact a well-defined form on $N$.  By Darboux's
Theorem there exist functions $x_1, y_1, p_1, q_1$ on $N$ such that
\[ d\tilde{\theta}_1 = -dp_1 \w dx_1 - dq_1 \w dy_1. \]
Then by Pfaff's Theorem there exists a function $z_1$ on $M_1$ such 
that
\[ \tilde{\theta}_1 = dz_1 - p_1\, dx_1 - q_1\, dy_1. \]
A similar argument shows that there exist functions $x_2, y_2, p_2, 
q_2$ on $N$ and $z_2, \mu$ on $M_2$ such that $\tilde{\theta}_2 = 
\mu\, \theta_2$ has the form
\[ \tilde{\theta}_2 = dz_2 - p_2\, dx_2 - q_2\, dy_2. \]

The ideal $\bar{\calI} = \{\om^1 \w \om^2,\ \om^3 \w \om^4\}$ on $N$ is 
now spanned by the forms 
\[ \{dp_1 \w dx_1 + dq_1 \w 
dy_1, \ dp_2 \w dx_2 + dq_2 \w dy_2\}, \]
and the ideals $\calI_1, \calI_2$ are {\em integrable extensions} of 
$\bar{\calI}$.  (For a definition and discussion of integrable 
extensions, see \cite{BG95}.)
The equations defining the B\"acklund transformation are simply those 
defining the change of coordinates on $N$:
\begin{align}
x_2 & = x_2(x_1, y_1, p_1, q_1) \notag \\
y_2 & = y_2(x_1, y_1, p_1, q_1) \label{Bzeroback} \\
p_2 & = p_2(x_1, y_1, p_1, q_1) \notag \\
q_2 & = q_2(x_1, y_1, p_1, q_1), \notag
\end{align}
or, in PDE notation,
\begin{align*}
\bar{x} & = \bar{x}(x, y, z_x, z_y) \\
\bar{y} & = \bar{x}(x, y, z_x, z_y) \\
\bar{z}_{\bar{x}} & = \bar{z}_{\bar{x}}(x, y, z_x, z_y) \\
\bar{z}_{\bar{y}} & = \bar{z}_{\bar{y}}(x, y, z_x, z_y).
\end{align*}
Note that if $z(x,y)$ is a known solution of the PDE corresponding to 
the ideal $(M_1, \calI_1)$, the corresponding solution $\bar{z}(\bar{x}, 
\bar{y})$ of the PDE corresponding to 
the ideal $(M_2, \calI_2)$ can be constructed by quadrature.

We have proved the following theorem:

\begin{theorem}\label{case3athm}
Let $\calB \subset M_1 \times M_2$ be a homogeneous B\"acklund 
transformation with the vectors $[C_1 \ \ C_2], \ [C_3 \ \ C_4]$ both 
nonzero and $B_1 = B_2 = B_3 = B_4 = 0.$  Then $\calB$ arises in the 
following way: let $\{x_1, y_1, p_1, q_1\}, \{x_2, y_2, p_2, q_2\}$ 
be two sets of local coordinates on a 4-manifold $N$ such that the 
2-forms
\[ \{dp_1 \w dx_1 + dq_1 \w dy_1, \ dp_2 \w dx_2 + dq_2 \w dy_2\} \]
span a hyperbolic pencil (i.e., there exist two distinct linear 
combinations of these 2-forms which are decomposable) at each point of 
$N$.  Let 
\begin{align*}
M_1 & = N \times \mathbb{R} \text{\ with coordinate $z_1$ on the 
$\mathbb{R}$ factor} \\
M_2 & = N \times \mathbb{R} \text{\ with coordinate $z_2$ on the 
$\mathbb{R}$ factor.}
\end{align*}
Let $\calI_1$ be the ideal on $M_1$ generated by the forms
\begin{gather*}
\theta_1 = dz_1 - p_1\, dx_1 - q_1\, dy_1 \\
d\theta_1 = -dp_1 \w dx_1 - dq_1 \w dy_1 \\
\Upsilon_1 = dp_2 \w dx_2 + dq_2 \w dy_2
\end{gather*}
and let $\calI_2$ be the ideal on $M_2$ generated by the forms
\begin{gather*}
\theta_2 = dz_2 - p_2\, dx_2 - q_2\, dy_2 \\
d\theta_2 = -dp_2 \w dx_2 - dq_2 \w dy_2 \\
\Upsilon_2 = dp_1 \w dx_1 + dq_1 \w dy_1.
\end{gather*}
Then $\calB \subset M_1 \times M_2$ is defined by the equations
\eqref{Bzeroback}.
\end{theorem}

Zvyagin \cite{Z81} calls B\"acklund transformations with $B_1 = B_2 = 
B_3 = B_4 = 0$ {\em holonomic.}  It can be shown that even without the 
assumption of homogeneity, any holonomic B\"acklund transformation 
arises locally from a hyperbolic system 
\[ \bar{\calI} = \{\om^1 \w \om^2, \, \om^3 \w \om^4\} \] 
on a 4-manifold $N$ such that $(M_1, \calI_1)$ 
and $(M_2, \calI_2)$ are integrable extensions of $(N, \bar{\calI})$.
These transformations are generally of limited interest.

\section{Case 3B: $B_1 = B_3 = 0;\ B_2, B_4 \neq 0$}

First we compute that
\begin{alignat*}{2}
0 & \equiv d(d\om^1) \mod{\{\om^1, \om^2\} &&} \\
& \equiv B_2 \, \theta_1 \w \theta_2 \w (G_3\, \om^3 + G_4\, \om^4) 
\qquad & \Rightarrow\  & G_3 = G_4 = 0 \\
0 & \equiv d(d\om^3) \mod{\{\om^3, \om^4\} &&} \\
& \equiv B_4 \, \theta_1 \w \theta_2 \w (H_1\, \om^1 + H_2\, \om^2) 
\qquad & \Rightarrow\  & H_1 = H_2 = 0 .
\end{alignat*}
Next we observe that under a transformation of the form 
\eqref{case3G1freedom}, we have
\begin{align*}
\tilde{B}_{2} & = b_{22} B_2 \\
\tilde{B}_{4} & = a_{22} B_4,
\end{align*} 
so we can choose a coframing with $B_2 = B_4 = 1$.  Such a coframing 
will be called {\em 2-adapted}; any two 2-adapted coframings differ by 
a transformation of the form
\begin{equation}
 \begin{bmatrix} \tilde{\theta}_1\\[0.1in] \tilde{\theta}_2\\[0.1in] 
\tilde{\om}^1\\[0.1in] \tilde{\om}^2\\[0.1in] \tilde{\om}^3\\[0.1in] 
\tilde{\om}^4 
\end{bmatrix} =  
\begin{bmatrix} 1 & 0 & 
0 & 0 & 0 & 0 \\[0.1in] 0 &  1 & 0 & 0 & 0 & 0 
\\[0.1in] 0 & 0 & 1 & 0 & 0 & 0 \\[0.1in] 0 & 0 & a_{21} & 
1 & 0 & 0 
\\[0.1in] 0 & 0 & 0 & 0 & 1 & 0 \\[0.1in] 0 & 0 & 0 & 
0 & b_{21} & 1
\end{bmatrix}^{-1} \begin{bmatrix} \theta_1\\[0.1in] \theta_2\\[0.1in] 
\om^1\\[0.1in] 
\om^2\\[0.1in] \om^3\\[0.1in] \om^4 \end{bmatrix}.  \label{case3BG2freedom}
\end{equation}
For a 2-adapted coframing, computing
\begin{align*}
d(d\om^2) \equiv 0 \mod{\{\om^1, \om^2\}} \\
d(d\om^4) \equiv 0 \mod{\{\om^3, \om^4\}}
\end{align*}
shows that
\begin{align*}
\alp_4 & = L_1\, \theta_1 + L_2\, \theta_2 - (A_1 + J_1 + 1)\, 
\om^1 - (J_2 + K_2)\, \om^2 + M_3\, \om^3 + M_4\, \om^4 \\
\bet_4 & = L_3\, \theta_1 + L_4\, \theta_2 + M_1\, \om^1 + M_2\, 
\om^2 - (A_2 + K_3 + 1)\, \om^3 - (J_4 + K_4)\, \om^4 
\end{align*}
for some functions $L_i, M_i$ on $\calB$.  

It is straightforward to show that under a transformation of 
the form \eqref{case3BG2freedom}, the functions $A_1, A_2, J_2, J_4, 
K_2, K_4, L_1, L_2, L_3, L_4, M_2, M_4$ remain unchanged.  By our 
assumption of homogeneity, they must therefore be constants.  
Moreover, 
\begin{align*}
\tilde{J}_1 & = J_1 + a_{21} J_2 \\
\tilde{K}_3 & = K_3 + b_{21} K_4 \\
\tilde{M}_1 & = M_1 + a_{21} M_2 \\
\tilde{M}_3 & = M_3 + b_{21} M_4 .
\end{align*}
Now we compute:
\begin{alignat*}{3}
0 & \equiv && d(d\theta_1) \mod{\{\theta_1, \theta_2, \om^3\}} && \\
& \equiv && A_1 (M_4 + 2 K_4)\, \om^1 \w \om^2 \w \om^4 & 
\Rightarrow &\  M_4 = -2K_4 \\
0 & \equiv && d(d\theta_2) \mod{\{\theta_1, \theta_2, \om^1\}} && \\
& \equiv && A_2 (M_2 + 2 J_2)\, \om^2 \w \om^3 \w \om^4 & 
\Rightarrow &\ M_2 = -2J_2 \\
0 & \equiv && d(d\theta_1) \mod{\{\theta_1, \om^3, \om^4\}} && \\
& \equiv &&A_1 (L_2 - L_4 - A_1 + 1)\, \theta_2 \w \om^1 \w \om^1  & 
\Rightarrow &\ L_4 = L_2 - A_1 + 1 \\
0 & \equiv && d(d\theta_2) \mod{\{\theta_2, \om^1, \om^2\}} && \\
& \equiv &&A_2 (L_3 - L_1 + A_2 - 1)\, \theta_1 \w \om^3 \w \om^4  & 
\Rightarrow &\ L_3 = L_1 - A_2 + 1 \\
0 & \equiv && d(d\theta_1) \mod{\{\theta_1, \theta_2, \om^4\}} && \\
& \equiv &&(A_1 M_3 + 2 A_1 K_3 + A_1 A_2 + A_1 - 1) &\, \om^1 \w \om^2 \w 
\om^3 &  \\
&&\ && \Rightarrow \ M_3 = & -2K_3 - A_2 - 1 - \frac{1}{A_1} \\
0 & \equiv && d(d\theta_2) \mod{\{\theta_1, \theta_2, \om^2\}} && \\
& \equiv && (A_2 M_1 + 2 A_2 J_1 + A_1 A_2 + A_2 - 1) & \om^1 \w \om^3 \w 
\om^4 &  \\
&&\ && \Rightarrow \ M_1 = & -2J_1 - A_1 - 1 - \frac{1}{A_2}  \\
0 & \equiv && d(d\theta_1) \mod{\{\om^1, \om^3, \om^4\}} && \\
& \equiv && J_2\, \theta_1 \w \theta_2 \w \om^2  & 
\Rightarrow &\ J_2 = 0 \\
0 & \equiv && d(d\theta_2) \mod{\{\om^1, \om^2, \om^3\}} && \\
& \equiv && K_4\, \theta_1 \w \theta_2 \w \om^4  & 
\Rightarrow &\ K_4 = 0.
\end{alignat*}
Therefore the functions $J_1, K_3$ in fact remain unchanged 
under a transformation of the form \eqref{case3BG2freedom}, and so 
they must be constants as well.  Next we compute:
\begin{alignat*}{2}
0 & \equiv d(d\om^2) \mod{\{\om^1, \om^3, \om^4\}} && \\
& \equiv (-L_1 - A_2 L_2)\, \theta_1 \w \theta_2 \w \om^2 \qquad & 
\Rightarrow\  & L_1 = -A_2 L_2 \\
0 & \equiv d(d\om^1) \mod{\{\om^2, \om^3, \om^4\}} && \\
& \equiv -K_2\, \theta_1 \w \theta_2 \w \om^1 \qquad & 
\Rightarrow\  & K_2 = 0 \\
0 & \equiv d(d\theta_1) \mod{\{\theta_2, \om^3, \om^4\}} && \\
& \equiv (A_1 A_2 - 1)(L_2 + 1)\, \theta_1 \w \om^1 \w \om^2 \qquad & 
\Rightarrow\  & L_2 = -1 \\
0 & \equiv d(d\om^3) \mod{\{\om^1, \om^2, \om^4\}} && \\
& \equiv -J_4\, \theta_1 \w \theta_2 \w \om^3 \qquad & 
\Rightarrow\  & J_4 = 0 \\
0 & \equiv d(d\theta_1) \mod{\{\om^2, \om^3, \om^4\}} && \\
& \equiv -(J_1 + A_1 + 1)\, \theta_1 \w \theta_2 \w \om^1 \qquad & 
\Rightarrow\  & J_1 = -A_1 - 1 \\
0 & \equiv d(d\theta_1) \mod{\{\om^1, \om^2, \om^4\}} && \\
& \equiv -A_1(K_3 + A_2 + 1)\, \theta_1 \w \theta_2 \w \om^3 \qquad & 
\Rightarrow\  & K_3 = -A_2 - 1 \\
0 & \equiv d(d\om^1) \mod{\{\theta_1, \om^2, \om^4\}} && \\
& \equiv (A_1 + 1)(A_2 + 1)\, \theta_2 \w \om^1 \w \om^3 \qquad & 
\Rightarrow\  & (A_1 + 1)(A_2 + 1) = 0.
\end{alignat*}
Without loss of generality, we can assume that $A_2 = -1$.  Then
\[ 0 = d(d\theta_2) = \frac{(A_1^2 - 1)}{A_1} \theta_2 \w \om^1 \w 
\om^3 \qquad \Rightarrow\  A_1^2 = 1. \]
Since $A_1 A_2 - 1 \neq 0$, we must have $A_1 = 1$.
The structure equations for a 2-adapted coframing may now be written 
as
\begin{align*}
d\theta_1 & = \theta_1 \w (\om^1 + \om^3) +  \om^1 \w \om^2 + \om^3 \w 
\om^4 \\
d\theta_2 & = -\theta_2 \w (\om^1 + \om^3) +  \om^1 \w \om^2 - \om^3 \w 
\om^4 \\
d\om^1 & = \om^1 \w (\theta_1 + \theta_2 - \om^3) \\
d\om^2 & = -\alp_3 \w \om^1 - \om^2 \w (\theta_1 + \theta_2 - 
\om^3) + \theta_1 \w \theta_2 + \om^3 \w \om^4 \\
d\om^3 & = \om^3 \w (\theta_2 - \theta_1 - \om^1) \\
d\om^4 & = -\bet_3 \w \om^3 - \om^4 \w (\theta_2 - \theta_1 - 
\om^1) + \theta_1 \w \theta_2 + \om^1 \w \om^2 
\end{align*}
for some 1-forms $\alp_3, \bet_3$ on $\calB$.

Now suppose that $\{\theta_1, \theta_2, \om^1, \om^2, \om^3, \om^4\}$ 
is any 2-adapted coframing.  Since 
\begin{align*}
d\om^1 & \equiv 0 \mod{\om^1} \\
d\om^3 & \equiv 0 \mod{\om^3}
\end{align*}
there exist functions $x, y, r_1, r_2$ on $\calB$ and nonzero constants 
$\lam_1, \lam_2$ such that
\[ \om^1 = \lam_1 e^{r_1}\, dx, \qquad \om^3 = \lam_2 e^{r_2}\, dy. \]
Since the systems $\{\theta_1, \om^1, \om^3\}$ and $\{\theta_2, \om^1, 
\om^3\}$ are completely integrable, there must exist functions $z_1, 
z_2, p_1, p_2, q_1, q_2, \rho_1, \rho_2$ on $\calB$, with $\rho_1, 
\rho_2$ nonvanishing, such that
\begin{align*}
\theta_1 & = \rho_1 (dz_1 - p_1\, dx - q_1\, dy)  \\
\theta_2 & = \rho_2 (dz_2 - p_2\, dx - q_2\, dy) . 
\end{align*}
Moreover, since 
\[ d\theta_1, d\theta_2 \equiv 0 \mod{\{\om^1, \om^3\}}, \]
$\rho_1$ must be a function of the variables $x, y, z_1$ alone and 
$\rho_2$ must be a function of the variables $x, y, z_2$ alone.
By making the contact transformation
\begin{gather*}
\tilde{x} = x \\
\tilde{y} = y \\
\tilde{z}_1 = -\tfrac{1}{2} \int_0^{z_1} \rho_1(x,y, \tau)\, 
d\tau \\
\tilde{z}_2 = -\tfrac{1}{2} \int_0^{z_2} \rho_2(x,y, \tau)\, 
d\tau  \\
\tilde{p}_1 = -\tfrac{1}{2} \int_0^{z_1} \frac{\partial \rho_1(x,y, 
\tau)}{\partial x}\, d\tau - \tfrac{1}{2} \rho_1(x,y,z_1)  p_1 \\ 
\tilde{p}_2 = -\tfrac{1}{2} \int_0^{z_2} \frac{\partial \rho_2(x,y, 
\tau)}{\partial x}\, d\tau - \tfrac{1}{2} \rho_2(x,y,z_2)  p_2 \\
\tilde{q}_1 = -\tfrac{1}{2} \int_0^{z_1} \frac{\partial \rho_1(x,y, 
\tau)}{\partial y}\, d\tau - \tfrac{1}{2} \rho_1(x,y,z_1)  q_1 \\ 
\tilde{q}_2 = -\tfrac{1}{2} \int_0^{z_2} \frac{\partial \rho_2(x,y, 
\tau)}{\partial y}\, d\tau - \tfrac{1}{2} \rho_2(x,y,z_2)  q_2
\end{gather*}
we can assume that $\rho_1 = \rho_2 = -\tfrac{1}{2}$.

Substituting the expressions given above for $\theta_1, \theta_2, \om^1, \om^3$ 
into the equations for $d\om^1, d\om^3$ yields
\[ \left. \begin{array}{l}
dr_1 \equiv \tfrac{1}{2} (dz_2 + dz_1) \\[0.1in]
dr_2 \equiv \tfrac{1}{2} (dz_2 - dz_1)
\end{array} \right\} \mod{\{dx, dy\}} . \]
Therefore we have
\begin{align*}
r_1 & = \tfrac{1}{2} (z_2 + z_1) + f(x, y) \\
r_2 & = \tfrac{1}{2} (z_2 - z_1) + g(x, y)
\end{align*}
for some functions $f, g$.  By making the contact transformation
\begin{gather*}
\tilde{x} = x \\
\tilde{y} = y \\
\tilde{z}_1 = z_1 + f(x,y) - g(x,y) \\
\tilde{z}_2 = z_2 + f(x,y) + g(x,y)  \\
\tilde{p}_1 = p_1 + \frac{\partial f}{\partial x} - \frac{\partial 
g}{\partial x} \\ 
\tilde{p}_2 = p_2 + \frac{\partial f}{\partial x} + \frac{\partial 
g}{\partial x} \\
\tilde{q}_1 = q_1 + \frac{\partial f}{\partial y} - \frac{\partial 
g}{\partial y} \\ 
\tilde{q}_2 = q_2 + \frac{\partial f}{\partial y} + \frac{\partial 
g}{\partial y}
\end{gather*}
we can assume that $f = g = 0$.  Now substituting into the equations 
for $d\om^1, d\om^3$ yields
\begin{align}
p_2 - p_1 & = 2\lam_1 e^{\frac{1}{2}(z_2 + z_1)} 
\label{liouvilleback} \\
q_2 + q_1 & = 2\lam_2 e^{\frac{1}{2}(z_2 - z_1)}. \notag
\end{align}
The equations for $d\theta_1, d\theta_2$ imply that
\begin{align*}
\om^2 & \equiv -\frac{1}{2\lam_1} e^{-\frac{1}{2}(z_2 + z_1)} dp_1 + 
\theta_1 \\
& \equiv -\frac{1}{2\lam_1} e^{-\frac{1}{2}(z_2 + z_1)} 
(dp_2 - 2\lam_1 \lam_2 e^{z_2}\, dy) - \theta_2 \mod{\om^1} \\
\om^4 & \equiv -\frac{1}{2\lam_2} e^{-\frac{1}{2}(z_2 - z_1)} dq_1 + 
\theta_1 \\
& \equiv \frac{1}{2\lam_2} e^{-\frac{1}{2}(z_2 - z_1)} (dq_2 -  
2\lam_1 \lam_2 e^{z_2}\, dx) + \theta_2 \mod{\om^2}.
\end{align*}
By scaling $x$ and $y$ if necessary, we can assume that
\[ \lam_1 = \frac{1}{2\lam_2} = \frac{\lam}{\sqrt{2}} \]
for some nonzero constant $\lam$.  Then equations 
\eqref{liouvilleback} become
\begin{align}
p_2 - p_1 & = \sqrt{2} \lam \, e^{\frac{1}{2}(z_2 + z_1)} 
\label{liouvilleback2} \\
q_2 + q_1 & = \frac{\sqrt{2}}{\lam}\,  e^{\frac{1}{2}(z_2 - z_1)}, 
\notag
\end{align}
or, in PDE notation,
\begin{align*}
\bar{z}_x - z_x & = \sqrt{2} \lam\,  e^{\frac{1}{2}(\bar{z} + z)} \\
\bar{z}_y + z_y & = \frac{\sqrt{2}}{\lam}\,  e^{\frac{1}{2}(\bar{z} - z)}. 
\end{align*}
This is the classical B\"acklund equation between the wave equation
\[ z_{xy} = 0 \]
and Liouville's equation
\[ \bar{z}_{xy} = e^{\bar{z}}. \]

We have proved the following theorem.

\begin{theorem}\label{case3bthm}
Let $\calB \subset M_1 \times M_2$ be a homogeneous B\"acklund 
transformation with the vectors $[C_1 \ \ C_2], \ [C_3 \ \ C_4], \ 
[B_1 \ \ B_2], \ [B_3 \ \ B_4]$ all nonzero and the pairs $[C_1 \ \ 
C_2], \ [B_1 \ \ B_2]$ and $[C_3 \ \ C_4], \ [B_3 \ \ B_4]$ both 
linearly dependent.  Then $\calB$ is locally contact equivalent to the 
transformation \eqref{liouvilleback2} between the wave equation 
\[ z_{xy} = 0 \]
and Liouville's equation
\[ \bar{z}_{xy} =  e^{\bar{z}}. \]
\end{theorem}

\section{Case 3C: Exactly one of $B_1, B_3$ is nonzero}

Without loss of generality, we can assume that $B_1 \neq 0, \ B_3 = 
0$.  Under a transformation of the form \eqref{case3G1freedom}, we have
\begin{align*}
\tilde{A}_1 & = \frac{b_{22}}{a_{22}} A_1 \\
\tilde{A}_2 & = \frac{a_{22}}{b_{22}} A_2 \\
\tilde{B}_1 & = (a_{22})^2 B_1 \\
\tilde{B}_2 & = -a_{21} a_{22} B_1 + b_{22} B_2 \\
\tilde{B}_4 & = a_{22} B_4.
\end{align*}
By Lemma \ref{Bveclemma}, the function $B_4$ must be nonzero, so we can 
choose a coframing with $B_2 = 0, \ A_1 = B_4 = 1.$  Such a coframing 
will be called {\em 2-adapted}.  By our homogeneity assumption, the 
functions $A_2, B_1$ are constant for any 2-adapted coframing.  
Moreover, any two 2-adapted coframings differ 
by a transformation of the form
\begin{equation}
 \begin{bmatrix} \tilde{\theta}_1\\[0.1in] \tilde{\theta}_2\\[0.1in] 
\tilde{\om}^1\\[0.1in] \tilde{\om}^2\\[0.1in] \tilde{\om}^3\\[0.1in] 
\tilde{\om}^4 
\end{bmatrix} =  
\begin{bmatrix} 1 & 0 & 
0 & 0 & 0 & 0 \\[0.1in] 0 &  1 & 0 & 0 & 0 & 0 
\\[0.1in] 0 & 0 & 1 & 0 & 0 & 0 \\[0.1in] 0 & 0 & 0 & 
1 & 0 & 0 
\\[0.1in] 0 & 0 & 0 & 0 & 1 & 0 \\[0.1in] 0 & 0 & 0 & 
0 & b_{21} & 1
\end{bmatrix}^{-1} \begin{bmatrix} \theta_1\\[0.1in] \theta_2\\[0.1in] 
\om^1\\[0.1in] 
\om^2\\[0.1in] \om^3\\[0.1in] \om^4 \end{bmatrix}.  \label{case3CG2freedom}
\end{equation}
For a 2-adapted coframing, computing
\[ d(d\om^2) \equiv 0 \mod{\{\om^1, \om^2\}} \]
shows that
\[ \alp_3 = L_1\, \theta_1 + L_2\, \theta_2 + M_1\, \om^1 + M_2\, 
\om^2 + \frac{1}{B_1} \, \om^3 \]
for some functions $L_i, M_i$ on $\calB$.
Using some of the remaining ambiguity \eqref{G0connfreedom}, we can assume 
that $M_1 = 0$.
Computing 
\begin{align*}
d(d\om^1) & \equiv 0 \mod{\{\om^1, \om^2\}} \\
d(d\om^4) & \equiv 0 \mod{\{\om^3, \om^4\}}
\end{align*}
shows that
\[ \alp_4 = P_1\, \theta_1 + P_2\, \theta_2 - (J_1 + 1)\, \om^1 + 
(B_1 - J_2 - K_2)\, \om^2 - \tfrac{1}{2}(A_2 + K_3)\, \om^3  -
\tfrac{1}{2} K_4\, \om^4 \]
for some functions $P_1, P_2$ on $\calB$.  Computing
\begin{align*}
d(d\theta_1) & \equiv 0 \mod{\theta_1} \\
d(d\theta_2) & \equiv 0 \mod{\theta_2}
\end{align*}
shows that
\begin{multline*} \bet_4 = (P_1 + 1)\, \theta_1 + (P_2 - 1)\, \theta_2 + 
\frac{(1 - A_2 - 2 A_2 J_1)}{A_2}\, \om^1 \\ + (B_1 - 2 J_2)\, \om^2 + 
\frac{(K_3 - A_2 - 2)}{2}\, \om^3 + \frac{(K_4 - 2 J_4)}{2}\, \om^4. 
\end{multline*}
It is straightforward to show that under a transformation of 
the form \eqref{case3BG2freedom}, the functions $G_4, J_1, J_2, 
J_4, K_2, K_4, L_1, L_2, M_2, P_1, P_2$ remain unchanged.  By our 
assumption of homogeneity, they must therefore be constants.  
Now we compute
\begin{alignat*}{2}
0 & \equiv d(d\om^3) \mod{\{\om^3, \om^4\}} && \\
& \equiv (H_1\, \om^1 + H_2\, \om^2) \w \theta_1 \w \theta_2) \qquad 
& \Rightarrow\  & H_1 = H_2 = 0 \\
0 & \equiv d(d\theta_1) \mod{\{\om^1, \om^2, \om^3\}} && \\
& \equiv \tfrac{1}{2} K_4\, \theta_1 \w \theta_2 \w \om^4 \qquad & 
\Rightarrow\  & K_4 = 0 \\
0 & \equiv d(d\theta_1) \mod{\{\om^1, \om^3, \om^4\}} && \\
& \equiv B_1(J_1 + 1)\, \theta_1 \w \theta_2 \w \om^2 \qquad & 
\Rightarrow\  & J_1 = -1 \\
0 & \equiv d(d\theta_1) \mod{\{\om^1, \om^3\}} && \\
& \equiv -G_4\, \theta_1 \w \om^2 \w \om^4 \qquad & 
\Rightarrow\  & G_4 = 0 \\
0 & \equiv d(d\theta_1) \mod{\{\om^2, \om^3\}} && \\
& \equiv \frac{(A_2 L_2 (B_1 - J_2) - P_2)}{A_2}\, \theta_1 \w \theta_2 
\w \om^1 \qquad & 
\Rightarrow\  & P_2 = A_2 L_2 (B_1 - J_2)  \\
0 & \equiv d(d\theta_2) \mod{\{\om^2, \om^3\}} && \\
& \equiv (-L_1 (B_1 - J_2) + P_1)\, \theta_1 \w \theta_2 \w \om^1 \qquad & 
\Rightarrow\  & P_1 = L_1 (B_1 - J_2).
\end{alignat*}
Next we compute
\begin{align*}
0 & \equiv d(d\theta_1) \mod{\om^3}  \\
& \equiv [(L_1 + A_2 L_2 + M_2)(J_2 - B_1) - K_2]\, \theta_1 \w \om^1 
\w \om^2 \\
0 & \equiv d(d\theta_2) \mod{\om^3} \\
& \equiv \left[(L_1 + A_2 L_2 + M_2)(J_2 - B_1) - \frac{K_2}{A_2}\right]\, 
\theta_2 \w \om^1 \w \om^2. 
\end{align*}
Since $A_2 \neq 1$, these equations imply that $K_2 = 0$.  Now
\begin{align*}
0 & \equiv d(d\om^4) \mod{\{\theta_2, \om^2, \om^3\}}  \\
& \equiv -(L_1 J_2 + 1)\, \theta_1 \w \om^1 \w \om^4.
\end{align*}
Therefore  $L_1, J_2$ are both nonzero, and 
\[ J_2 = -\frac{1}{L_1}. \]
Next we compute
\begin{alignat*}{2}
0 & \equiv d(d\om^4) \mod{\{\theta_2, \om^1, \om^3\}} && \\
& \equiv B_1 (A_2 + 1)\, \theta_1 \w \om^2 \w \om^4 \qquad & 
\Rightarrow\  & A_2 = -1 \\
0 & \equiv d(d\om^4) \mod{\{\theta_1, \om^2, \om^3\}} && \\
& \equiv \frac{(L_2 - L_1)}{L_1}\, \theta_2 \w \om^1 \w \om^4 \qquad & 
\Rightarrow\  & L_2 = L_1\\
0 & \equiv d(d\om^2) \mod{\om^3} && \\
& \equiv B_1 M_2\, \theta_1 \w \theta_2 \w \om^2  \qquad & 
\Rightarrow\  & M_2 = 0.
\end{alignat*}
Now it is straightforward to check that under a transformation of 
the form \eqref{case3BG2freedom}, the functions $G_3, K_3$ remain 
unchanged; therefore they must be constants.  Continuing, we have
\begin{alignat*}{2}
0 & \equiv d(d\om^3) \mod{\{\theta_1, \om^4\}} && \\
& \equiv -J_4\, \om^1 \w \om^2 \w \om^3 \qquad & 
\Rightarrow\  & J_4 = 0 \\
0 & \equiv d(d\om^1) \mod{\{\theta_1, \theta_2\}} && \\
& \equiv \frac{(1 - K_3)}{L_1}\, \om^1 \w \om^2 \w \om^3 \qquad & 
\Rightarrow\  & K_3 = 1\\
0 & \equiv d(d\theta_1) \mod{\{\theta_2, \om^2, \om^4\}} && \\
& \equiv -\frac{(L_1 B_1 + 1)}{L_1 B_1}\, \theta_1 \w \om^1 \w \om^3 \qquad & 
\Rightarrow\  & L_1 = -\frac{1}{B_1}\\
0 & = d(d\theta_1)  && \\
& = (B_1 - G_3)\, \theta_1 \w \om^2 \w \om^3 \qquad & 
\Rightarrow\  & G_3 = B_1.
\end{alignat*}

In summary, we have now shown that the structure equations 
of a 2-adapted coframing take the form
\begin{align}
d\theta_1 & = \theta_1 \w (\om^1 + \om^3) +  \om^1 \w \om^2 + \om^3 \w 
\om^4 \notag \\
d\theta_2 & = -\theta_2 \w (\om^1 + \om^3) +  \om^1 \w \om^2 - \om^3 \w 
\om^4 \notag \\
d\om^1 & = B_1 (\theta_1 \w \theta_2 + \theta_1 \w \om^2 + \theta_2 
\w \om^2 + \om^2 \w \om^3) \label{case3Cstreqs} \\
d\om^2 & = \frac{1}{B_1}\, (\theta_1 + \theta_2 - \om^3) \w \om^1 + 
\om^3 \w \om^4 \notag \\
d\om^3 & = (\theta_1 - \theta_2 - B_1\, \om^2) \w \om^3 \notag \\
d\om^4 & = -\bet_3 \w \om^3 + \theta_1 \w \theta_2 - \theta_1 \w \om^4 + 
\theta_2 \w \om^4 + B_1\, \om^2 \w \om^4 + \om^1 \w \om^2 \notag
\end{align}
for some 1-form $\bet_3$ on $\calB$.

\begin{lemma}
We can choose a 2-adapted coframing with $\bet_3 = 0$.
\end{lemma}

\begin{pf}
We will make liberal use of the fact that $\bet_3$ is only 
well-defined modulo $\om^3$, so we can add multiples of $\om^3$ 
to $\bet_3$ at will.

Suppose that $\bet_3 \neq 0$.  Differentiating the last equation in 
\eqref{case3Cstreqs} yields
\[ 0 = d(d\om^4) = [-d\bet_3 + -2(\theta_1 - \theta_2 - B_1\, 
\om^2) \w \bet_3] \w \om^3; \]
therefore
\begin{equation}
d\bet_3 \equiv -2(\theta_1 - \theta_2 - B_1\, \om^2) \w \bet_3 
\mod{\om^3}. \label{dbet3eq}
\end{equation}
By the Frobenius Theorem, there exist functions $x, y, \lambda, \mu, 
\nu$ on $\calB$ such that
\begin{align*}
\om^3 & = e^{\lambda}\, dx \\
\bet_3 & = e^{\mu}\, dy + \nu\, dx.
\end{align*}
In fact, because of the $\om^3$-ambiguity in $\bet_3$, we can assume 
that $\nu = 0$.

From the equation for $d\om^3$ in \eqref{case3Cstreqs} and equation 
\eqref{dbet3eq} we have
\begin{align*}
[d\lambda - (\theta_1 - \theta_2 - B_1\, \om^2)] \w \om^3 & = 0 \\
[d\mu + 2 (\theta_1 - \theta_2 - B_1\, \om^2)]  \w \bet_3 & \equiv 0 
\mod{\om^3}.
\end{align*}
It follows that
\begin{align*}
d\lambda & = \theta_1 - \theta_2 - B_1\, \om^2 + r_1\, dx \\
d\mu & = -2 (\theta_1 - \theta_2 - B_1\, \om^2) + r_2\, dx + r_3\, dy
\end{align*}
for some functions $r_1, r_2, r_3$ on $\calB$.  Therefore
$d(\mu + 2\lambda)$ is a linear combination of $dx$ and $dy$.  This 
implies that there exists a function $f(x,y)$ such that
\[ \mu + 2\lambda = f(x,y). \]
Thus we have
\[ \bet_3 = e^{-2\lambda} e^{f(x,y)} dy, \]
and by replacing the function $y$ by the function $\int e^{f(x,y)} 
dy$ (and adding multiples of $\om^3$ to $\bet_3$ if necessary) we 
can assume that
\[ \bet_3 = e^{-2\lambda} dy. \]

From the structure equation for $d\om^4$ in \eqref{case3Cstreqs}, it 
follows that under a change of 2-adapted coframing of the form
\[ \tilde{\om}^4 = \om^4 - b_{21}\, \om^3 \]
we have
\begin{align*}
\tilde{\bet}_3 & = \bet_3 + db_{21} + 2b_{21}\, d\lambda \\
& = e^{-2\lambda} dy + db_{21} + 2b_{21}\, d\lambda . 
\end{align*}
Taking $b_{21} = -ye^{-2\lambda}$ yields $\tilde{\bet}_3 = 0$, as desired.
\end{pf}

For a 2-adapted coframing as in the lemma, the structure equations 
\eqref{case3Cstreqs} take the form
\begin{align*}
d\theta_1 & = \theta_1 \w (\om^1 + \om^3) +  \om^1 \w \om^2 + \om^3 \w 
\om^4 \\
d\theta_2 & = -\theta_2 \w (\om^1 + \om^3) +  \om^1 \w \om^2 - \om^3 \w 
\om^4 \\
d\om^1 & = B_1 (\theta_1 \w \theta_2 + \theta_1 \w \om^2 + \theta_2 
\w \om^2 + \om^2 \w \om^3) \\
d\om^2 & = \frac{1}{B_1}\, (\theta_1 + \theta_2 - \om^3) \w \om^1 + 
\om^3 \w \om^4 \\
d\om^3 & = (\theta_1 - \theta_2 - B_1\, \om^2) \w \om^3 \\
d\om^4 & = \theta_1 \w \theta_2 - \theta_1 \w \om^4 + \theta_2 \w 
\om^4 + B_1\, \om^2 \w \om^4 + \om^1 \w \om^2.
\end{align*}
The intepretation of these equations requires some 
preliminaries regarding the geometry of frame bundles, which we will 
postpone until after the next section.

\section{Case 3D: $B_1, B_3 \neq 0$}

Under a transformation of the form \eqref{case3G1freedom}, we have
\begin{align*}
\tilde{A}_1 & = \frac{b_{22}}{a_{22}} A_1 \\
\tilde{A}_2 & = \frac{a_{22}}{b_{22}} A_2 \\
\tilde{B}_1 & = (a_{22})^2 B_1 \\
\tilde{B}_2 & = -a_{21} a_{22} B_1 + b_{22} B_2 \\
\tilde{B}_3 & = (b_{22})^2 B_3 \\
\tilde{B}_4 & = -b_{21} b_{22} B_3 + a_{22} B_4.
\end{align*}
Since $B_1, B_3 \neq 0$, we can choose a coframing for which $B_2 = 
B_4 = 0$ and $B_1,\, B_3$ are constants.  Such a coframing 
will be called {\em 2-adapted.}  A 2-adapted coframing 
is uniquely determined by the constants $B_1$ and $B_3$, and by the 
homogeneity assumption, all the other torsion functions are constants as well.
For now we will not specify the values of $B_1$ and $B_3$; rather 
we will use this ambiguity to specify the values of other torsion 
coefficients in what follows.  

For a 2-adapted coframing, computing 
\begin{align*}
d(d\om^2) & \equiv 0 \mod{\{\om^1, \om^2\}} \\
d(d\om^4) & \equiv 0 \mod{\{\om^3, \om^4\}}
\end{align*}
shows that
\begin{align*}
\alp_3 & = L_1\, \theta_1 + L_2\, \theta_2 + M_1\, \om^1 + M_2\, 
\om^2 - \frac{B_3}{B_1}\, \om^4 \\
\bet_3 & = L_3\, \theta_1 + L_4\, \theta_2 - \frac{B_1}{B_3}\, \om^2 + 
M_3\, \om^3 - M_4\, \om^4
\end{align*}
for some functions $L_i, M_i$ on $\calB$.  Using the remaining 
ambiguity \eqref{G0connfreedom}, we can assume that $M_1 = M_3 = 0$.  
Computing 
\begin{align*}
d(d\om^1) & \equiv 0 \mod{\{\om^1, \om^2\}} \\
d(d\om^3) & \equiv 0 \mod{\{\om^3, \om^4\}}
\end{align*}
shows that
\begin{align*} 
\alp_4 & = P_1\, \theta_1 + P_2\, \theta_2 + Q_1\, \om^1 + 
Q_2\, \om^2 - \tfrac{1}{2}(A_2 + K_3)\, \om^3  -
\tfrac{1}{2} K_4\, \om^4 \\
\bet_4 & = P_3\, \theta_1 + P_4\, \theta_2 - \tfrac{1}{2}(A_1 + J_1)\, 
\om^1 - \tfrac{1}{2} J_2\, \om^2 + Q_3\, \om^3 + Q_4\, \om^4
\end{align*}
for some functions $P_i, Q_i$ on $\calB$.  Now we compute
\begin{alignat*}{2}
0 & \equiv d(d\theta_1) \mod{\{\theta_1, \om^3, \om^4\}} && \\
& \equiv  A_1 (P_2 - P_4)\, \theta_2 \w \om^1 \w \om^2 \qquad 
& \Rightarrow\  & P_4 = P_2 \\
0 & \equiv d(d\theta_2) \mod{\{\theta_2, \om^1, \om^2\}} && \\
& \equiv  A_2 (P_3 - P_1)\, \theta_1 \w \om^3 \w \om^4 \qquad 
& \Rightarrow\  & P_3 = P_1 \\
0 & \equiv d(d\theta_1) \mod{\{\theta_1, \theta_2, \om^3\}} && \\
& \equiv  A_1 (\tfrac{1}{2} K_4 - J_4 - Q_4)\, \om^1 \w \om^2 \w 
\om^4 \qquad & \Rightarrow\  & Q_4 = \tfrac{1}{2} K_4 - J_4 \\
0 & \equiv d(d\theta_2) \mod{\{\theta_1, \theta_2, \om^1\}} && \\
& \equiv  A_2 (\tfrac{1}{2} J_2 - K_2 - Q_2)\, \om^2 \w \om^3 \w 
\om^4 \qquad & \Rightarrow\  & Q_2 = \tfrac{1}{2} J_2 - K_2 \\
0 & \equiv d(d\theta_1) \mod{\{\om^1, \om^3, \om^4\}} && \\
& \equiv  -B_1 Q_1 \, \theta_1 \w \theta_2 \w 
\om^2 \qquad & \Rightarrow\  & Q_1 = 0\\
0 & \equiv d(d\theta_2) \mod{\{\om^1, \om^2, \om^3\}} && \\
& \equiv  -B_3 Q_3 \, \theta_1 \w \theta_2 \w 
\om^4 \qquad & \Rightarrow\  & Q_3 = 0\\
0 & \equiv d(d\theta_1) \mod{\{\om^1, \om^2, \om^4\}} && \\
& \equiv  (\tfrac{1}{2} K_4 L_4 - A_2 P_2)\, \theta_1 \w \theta_2 \w 
\om^3 \qquad & \Rightarrow\  & P_2 = \frac{K_4 L_4}{2 A_2} \\
0 & \equiv d(d\theta_2) \mod{\{\om^2, \om^3, \om^4\}} && \\
& \equiv  (A_1 P_1 - \tfrac{1}{2} J_2 L_1)\, \theta_1 \w \theta_2 \w 
\om^1 \qquad & \Rightarrow\  & P_1 = \frac{J_2 L_1}{2 A_1} \\
0 & \equiv d(d\theta_1) \mod{\theta_1} && \\
& \equiv  (\tfrac{1}{2} A_1 K_3 - \tfrac{1}{2} A_1 A_2 - 1)\, \om^1 \w 
\om^2 \w \om^3 \qquad & \Rightarrow\  & K_3 = \frac{A_1 A_2 + 2}{A_1} \\
0 & \equiv d(d\theta_2) \mod{\theta_2} && \\
& \equiv  (\tfrac{1}{2} A_2 J_1 - \tfrac{1}{2} A_1 A_2 - 1)\, \om^1 \w 
\om^3 \w \om^4 \qquad & \Rightarrow\  & J_1 = \frac{A_1 A_2 + 2}{A_2} \\
0 & \equiv d(d\theta_1) \mod{\{\theta_2, \om^2, \om^4\}} && \\
& \equiv \frac{(1 - (A_1 A_2)^2)}{A_1 A_2}\, \theta_1 \w 
\om^1 \w \om^3 \qquad & \Rightarrow\  & (A_1 A_2)^2 = 1. 
\end{alignat*}
Since we require $A_1 A_2 \neq 1$, we must have $A_1 A_2 = -1$, and 
therefore $A_2 = -\frac{1}{A_1}$.  Continuing, we have
\begin{alignat*}{2}
0 & \equiv d(d\theta_1) \mod{\{\theta_2, \om^1, \om^3\}} && \\
& \equiv  \frac{(H_2 - A_1^2 G_4)}{A_1}\, \theta_1 \w \om^2 \w \om^4 \qquad 
& \Rightarrow\  & H_2 = A_1^2 G_4 \\
0 & \equiv d(d\om^1) \mod{\{\om^2, \om^3, \om^4\}} && \\
& \equiv  \frac{B_1 (L_2 - A_1 L_1)}{A_1}\, \theta_1 \w \theta_2 \w \om^1 \qquad 
& \Rightarrow\  & L_2 = A_1 L_1 \\
0 & \equiv d(d\om^3) \mod{\{\om^1, \om^2, \om^4\}} && \\
& \equiv  -B_3(L_4 + A_1 L_3) \, \theta_1 \w \theta_2 \w \om^3 \qquad 
& \Rightarrow\  & L_4 = -A_1 L_3 \\
0 & \equiv d(d\om^1) \mod{\{\theta_1, \theta_2, \om^4\}} && \\
& \equiv  \frac{(B_3 G_3 M_2 - B_1 K_4)}{B_3}\, \om^1 \w \om^2 \w \om^3 \qquad 
& \Rightarrow\  & K_4 = \frac{B_3 G_3 M_2}{B_1} \\
0 & \equiv d(d\om^3) \mod{\{\theta_1, \theta_2, \om^2\}} && \\
& \equiv  \frac{(B_1 H_1 M_4 - B_3 J_2)}{B_1}\, \om^1 \w \om^3 \w \om^4 \qquad 
& \Rightarrow\  & J_2 = \frac{B_1 H_1 M_4}{B_3} \\
0 & \equiv d(d\om^2) \mod{\{\om^1, \om^3, \om^4\}} && \\
& \equiv B_1 M_2 \, \theta_1 \w \theta_2 \w \om^2 \qquad 
& \Rightarrow\  & M_2 = 0 \\
0 & \equiv d(d\om^4) \mod{\{\om^1, \om^2, \om^3\}} && \\
& \equiv B_3 M_4 \, \theta_1 \w \theta_2 \w \om^4 \qquad 
& \Rightarrow\  & M_4 = 0.
\end{alignat*}
Now we have
\begin{gather*}
0 = d(d\theta_1) = [-A_1 K_2\, \om^1 \w \om^2 + \frac{H_1}{A_1}\, \om^1 \w 
\om^4 - A_1 G_3\, \om^2 \w \om^3 - \frac{J_4}{A_1}\, \om^3 \w \om^4] 
\w \theta_1 \\
\Rightarrow\  G_3 = H_1 = J_4 = K_2 = 0
\end{gather*}
\begin{alignat*}{2}
0 & \equiv d(d\om^2) \mod{\{\theta_2, \om^2, \om^4\}} && \\
& \equiv  \frac{(B_1 L_1 + A_1 B_3 L_3)}{A_1 B_1}\, \theta_1 \w \om^1 \w 
\om^3 \qquad & \Rightarrow\  & L_1 = -\frac{A_1 B_3 L_3}{B_1} \\
0 & \equiv d(d\om^2) \mod{\{\theta_1, \om^1, \om^3\}} && \\
& \equiv  \frac{-B_3(B_1 + A_1^2 G_4 L_3)}{B_1}\, \theta_2 \w \om^2 \w 
\om^4 \qquad & \Rightarrow\  & B_1 + A_1^2 G_4 L_3 = 0. 
\end{alignat*}
Since $B_1 \neq 0$, this implies that $G_4, L_3 \neq 0$ and $L_3 = 
-\frac{B_1}{A_1^2 G_4}$.

The structure equations now take the form
\begin{align*}
d\theta_1 & = \theta_1 \w (A_1\, \om^1 + \frac{1}{A_1}\, \om^3) +  
A_1\, \om^1 \w \om^2 + \om^3 \w \om^4 \\
d\theta_2 & = -\theta_2 \w (A_1\, \om^1 + \frac{1}{A_1}\, \om^3) +  
\om^1 \w \om^2 - \frac{1}{A_1}\, \om^3 \w \om^4 \\
d\om^1 & = (\frac{B_1}{A_1}\, \theta_1 + B_1\, \theta_2 - 
G_4\, \om^4) \w \om^2 + B_1\, \theta_1 \w \theta_2 \\
d\om^2 & = (-\frac{B_3}{A_1 G_4}\, \theta_1 - \frac{B_3}{G_4}\, 
\theta_2 + \frac{B_3}{B_1}\, \om^4) \w \om^1 + \om^3 \w \om^4  \\
d\om^3 & = (-B_3\, \theta_1 + A_1 B_3\, \theta_2 - A_1^2 G_4\, \om^2) 
\w \om^4 + B_3\, \theta_1 \w \theta_2  \\
d\om^4 & = (\frac{B_1}{A_1^2 G_4}\, \theta_1 - \frac{B_1}{A_1 G_4}\, 
\theta_2 + \frac{B_1}{B_3}\, \om^2) \w \om^3 + \om^1 \w \om^2.
\end{align*}
By a transformation of the form \eqref{case3G1freedom} with 
\[ a_{22} = \frac{1}{\sqrt{|G_4|}}, \qquad b_{22} = \frac{1}{A_1 
\sqrt{|G_4|}} \]
we can arrange that $A_1 = 1, G_4 = \pm 1$.  Let $\eps = G_4 = \pm 
1$; then the structure equations take the form
\begin{align*}
d\theta_1 & = \theta_1 \w (\om^1 + \om^3) +  
\om^1 \w \om^2 + \om^3 \w \om^4 \\
d\theta_2 & = -\theta_2 \w (\om^1 + \om^3) +  
\om^1 \w \om^2 - \om^3 \w \om^4 \\
d\om^1 & = (B_1\, \theta_1 + B_1\, \theta_2 - 
\eps\, \om^4) \w \om^2 + B_1\, \theta_1 \w \theta_2 \\
d\om^2 & = (-\eps B_3\, \theta_1 - \eps B_3\, 
\theta_2 + \frac{B_3}{B_1}\, \om^4) \w \om^1 + \om^3 \w \om^4  \\
d\om^3 & = (-B_3\, \theta_1 + B_3\, \theta_2 - \eps\, \om^2) 
\w \om^4 + B_3\, \theta_1 \w \theta_2  \\
d\om^4 & = (\eps B_1\, \theta_1 - \eps B_1\, 
\theta_2 + \frac{B_1}{B_3}\, \om^2) \w \om^3 + \om^1 \w \om^2.
\end{align*}
The next section contains a discussion of frame bundles which will be 
necessary in order to interpret these structure equations and those of 
the previous section.

\section{Local geometry of surfaces in 3-dimensional Riemannian and 
Lorentzian space forms}

First we will discuss the familiar geometry of surfaces in 
3-dimensional Euclidean space; then we can examine what changes when 
the curvature and/or the signature of the underlying space form is allowed 
to vary.

Let $\bb{E}^3$ denote the vector space $\bb{R}^3$ with the 
Euclidean inner product
\[ \langle x, y \rangle = x^1 y^1 + x^2 y^2 + x^3 y^3. \]
An {\em orthonormal frame} at a point $x \in \bb{E}^3$ an orthonormal 
basis $\{e_1, e_2, e_3\}$ for the tangent space $T_x\bb{E}^3$.  The 
set of all orthonormal frames at all points of $\bb{E}^3$ is called 
the {\em frame bundle} of $\bb{E}^3$, denoted $\calF(\bb{E}^3)$; it is 
a principal bundle over 
$\bb{E}^3$ whose fiber over each point $x \in \bb{E}^3$ is naturally 
isomorphic to the Lie group $O(3)$ (or, if we require our frames to be 
positively oriented, $SO(3)$).

The frame bundle $\calF(\bb{E}^3)$ is in fact naturally isomorphic to 
the Lie group $E(3)$, the group of isometries of $\bb{E}^3$.  Recall 
that
\[ E(3) = \left\{ \begin{bmatrix} A & b \\ 0 & 1 \end{bmatrix} : A \in 
O(3), b \in \bb{E}^3 \right\}. \]
If we represent a vector 
$y \in \bb{E}^3$ by the 4-dimensional vector 
$\begin{bmatrix} y \\ 1 \end{bmatrix}$, then elements of $E(3)$ 
act on $y$ by matrix multiplication:
\[ \begin{bmatrix} A & b \\ 0 & 1 \end{bmatrix} \begin{bmatrix} y \\ 1 
\end{bmatrix} = \begin{bmatrix} Ay + b \\ 1 \end{bmatrix}. \]  
An orthonormal frame $\{e_1, e_2, e_3\}$ at $x \in \bb{E}^3$ may be 
regarded as an element of $E(3)$ by letting $A$ be the matrix whose 
columns are the vectors $e_1, e_2, e_3$ and letting $b$ be the vector $x$.

The vectors $x, e_1, e_2, e_3$ may all be thought of as 
$\bb{E}^3$-valued 
functions on $\calF(\bb{E}^3)$.  Thus their exterior derivatives $dx, de_i$ are
$T\bb{E}^3$-valued 1-forms on $\calF(\bb{E}^3)$.  Since $\{e_1, e_2, e_3 \}$ 
is a basis for the 
tangent space to $\bb{E}^3$ at each point, we can express $dx, de_i$ 
as linear combinations of $e_1, e_2, e_3$ whose 
coefficients are ordinary scalar-valued 1-forms on $\calF(\bb{E}^3)$.  
Hence we can define 1-forms $\eta^i,\, \eta^j_i, \ 1 \leq i,j \leq 3$, on 
$\calF(\bb{E}^3)$ by the equations
\begin{align}
dx & = \sum_{i=1}^3 e_i\, \eta^i \label{vecstruceuc}  \\ 
de_i & = \sum_{j=1}^3 e_j\, \eta^j_i. \notag
\end{align} 
The 1-forms $\eta^1, \eta^2, \eta^3$ are semi-basic for the natural 
projection $\pi: \calF(\bb{E}^3) \to \bb{E}^3$.  They have the 
property that if $\sigma: \bb{E}^3 \to \calF(\bb{E}^3)$ is a section 
of the frame bundle defined by
\[ \sigma(x) = (e_1(x), e_2(x), e_3(x)), \]
then the pullbacks $\sigma^*(\eta^i)$ are dual to the basis $\{e_1(x), 
e_2(x), e_3(x)\}$ for the tangent space $T_x\bb{E}^3$ at each point $x 
\in \bb{E}^3$.  Thus the forms 
$\{\sigma^*(\eta^1), \sigma^*(\eta^2), 
\sigma^*(\eta^3)\}$ are a basis for the 1-forms on $\bb{E}^3$.  For this 
reason, the $\eta^i$ are called the {\em dual forms} on $\calF(\bb{E}^3)$.
The $\eta^i_j$, on the other hand, form a basis for the 1-forms on each 
fiber of $\pi$.  
If $\sigma$ is a section as above, then the pullbacks $\sigma^*(\eta^i_j)$ 
are the Levi-Civita connection forms of the Euclidean metric on 
$\bb{E}^3$ for the frame defined by 
$\sigma$.  For this reason, the $\eta^i_j$ are called the {\em 
connection forms} on $\calF(\bb{E}^3)$.  Together, the forms $\{\eta^i, 
\eta^i_j\}$ form a basis for the left-invariant forms on the group 
$E(3)$, and hence for the Lie algebra $\mathfrak{e}(3)$.


Differentiating equations \eqref{vecstruceuc} shows that the forms 
$\eta^i,\, \eta^j_i$ satisfy the {\em structure equations}
\begin{align}
d\eta^i & = - \sum_{j=1}^3 \eta^i_j \w \eta^j \label{formstruceuc} \\
d\eta^i_j & = - \sum_{k=1}^3 \eta^i_k \w \eta^k_j. \notag
\end{align}
(These equations are equivalent to the structure equations for the Lie 
algebra $\mathfrak{e}(3)$.)
Differentiating the equations
\[ \langle e_i, e_j \rangle = 
\begin{cases} 
0 & i \neq j \\
1 & i = j 
\end{cases} \]
shows that $\eta^i_j = -\eta^j_i$; in particular, $\eta^i_i = 0$ for $i = 
1,2,3$.  

Now let $U \subset \bb{R}^2$ be open and let
$X:U \to \bb{E}^3$ be a regular surface.  An {\em adapted
orthonormal frame field} along $X$ is a choice, for each $x \in X$, of 
an orthonormal frame $\{e_1, e_2, e_3\}$ at $x$ such that the vectors 
$e_1, e_2$ form a basis for the tangent space $T_xX$ (and hence
$e_3$ is a unit normal vector to $X$ at $x$.)
If $\{e_1, e_2, e_3\}$ is an adapted 
orthonormal frame field, then any other adapted orthonormal frame 
field $\{\tilde{e}_1, \tilde{e}_2, \tilde{e}_3 \}$ has the form
\begin{align*}
\tilde{e}_1 & = \pm [(\cos \varphi)\, e_1 - (\sin \varphi)\, e_2] \\
\tilde{e}_2 & = \pm [(\sin \varphi)\, e_1 + (\cos \varphi)\, e_2] \\
\tilde{e}_3 & = \pm e_3 
\end{align*}
for some function $\varphi$ on $X$.  (The ambiguities of sign can be 
removed by specifying a choice of unit normal and requiring that the 
frame field be positively oriented.)

A choice of an adapted orthonormal frame field may be thought of as a 
lifting $\tilde{X}:U \to \calF(\bb{E}^3)$.  Now consider the pullbacks of the 
forms $\eta^i, \eta^i_j$ via $\tilde{X}$ to the surface $X$.  (The pullback 
notation will be omitted for simplicity.)
Since $e_1, e_2$ form a basis for $T_xX$ at each point $x \in X$, the 
1-form $dx = \sum e_i\, \eta^i$ must be a linear combination of $e_1$ and $e_2$; 
therefore, $\eta^3 = 0$. Moreover, the 1-forms $\eta^1, \eta^2$ are 
linearly independent and so form a basis for the 1-forms on $X$.  
Differentiating the equation $\eta^3 = 0$ yields
\[ 0 = d\eta^3 = -\eta^3_1 \w \eta^1 - \eta^3_2 \w \eta^2. \]
By Cartan's Lemma, there exist functions $h_{11}, h_{12}, h_{22}$ on 
$X$ such that
\begin{gather*}
\eta^3_1 = h_{11}\, \eta^1 + h_{12}\, \eta^2 \\
\eta^3_2 = h_{12}\, \eta^1 + h_{22}\, \eta^2.
\end{gather*}
The structure equations for the dual forms can now be written in the form
\begin{gather}
d\eta^1 = -\eta^1_2 \w \eta^2 \label{streqs} \\
d\eta^2 = \eta^1_2 \w \eta^1 \notag
\end{gather}
where $\eta^1_2$ is the Levi-Civita connection form for the induced metric on 
$X$.  
The first and second fundamental forms of $X$ are
\begin{gather*}
I = \langle dX, dX \rangle = (\eta^1)^2 + (\eta^2)^2 \\
II = \langle dX, de_3 \rangle = h_{11}\, (\eta^1)^2 + 2h_{12}\,\eta^1 
\, \eta^2 + h_{22}\, (\eta^2)^2.
\end{gather*}
The {\em Gauss curvature} $K$ of $X$ is defined to be the determinant of 
$II$, i.e.,
\[ K = h_{11}h_{22} - h_{12}^2, \]
and the structure equation for the Levi-Civita form $\eta^1_2$ can now 
be written in the form
\[ d\eta^1_2 = K\, \eta^1 \w \eta^2. \]
(This is called the {\em Gauss equation}.)
The {\em mean curvature} $H$ of $X$ is defined to be one-half of the trace 
of $II$ with respect to the metric defined by $I$, i.e., 
\[ H = \tfrac{1}{2}(h_{11} + h_{22}). \]
The quantities $K$ and $H$ are independent (up to the sign of $H$) of the choice 
of adapted orthonormal frame field on $X$.  Note that
\begin{gather*}
\eta^3_1 \w \eta^3_2 = (h_{11}h_{22} - h_{12}^2)\, \eta^1 
\w \eta^2  = K\, \eta^1 \w \eta^2 \\
\eta^3_1 \w \eta^2 + \eta^1 \w \eta^3_2 = (h_{11} + h_{22}) \,\eta^1 \w 
\eta^2 = 2H\, \eta^1 \w \eta^2. 
\end{gather*}

So for instance, let $X: U \to \bb{E}^3$ be any surface whose Gauss 
curvature $K$ satisfies $K \equiv -1$.  If $\tilde{X}:U \to 
\calF(\bb{E}^3)$ is any choice of adapted orthonormal coframing along 
$X$, then the image of $\tilde{X}$ is an integral manifold of the 
exterior differential system 
\[ \tilde{\calI} = \{\eta^3,\ d\eta^3,\ \eta^3_1 \w \eta^3_2 + \eta^1 \w 
\eta^2 \}  \]
on $\calF(\bb{E}^3)$.
But there is one further wrinkle to consider.  Generally our objects 
of interest are surfaces, and while the unit normal vector $e_3$ is 
determined up to sign by the surface, in general there is no canonical 
choice of basis $\{e_1, e_2\}$ for the tangent spaces $T_xX$.  
Rather than lifting $X$ to the entire frame bundle 
$\calF(\bb{E}^3)$, it is more natural to consider liftings of $X$ to 
the space $M$ of {\em contact elements} of $\bb{E}^3$.  This is the space 
of tangent planes to points of $\bb{E}^3$, and if we allow these 
planes to be oriented by a choice of unit normal vector, $M$ may 
be described as
\[ M = \{(x, e_3): x \in \bb{E}^3, e_3 \in T_x\bb{E}^3, \langle e_3, 
e_3 \rangle = 1 \}. \]
This is a 5-dimensional manifold, and it is naturally the quotient of 
$\calF(\bb{E}^3)$ by the circle action consisting of rotations between 
$e_1$ and $e_2$ at each point.

The 1-form $\eta^3$ is well-defined on $M$, and in fact it is a 
contact form on $M$.  The forms 
$\eta^1, \eta^2, \eta^3_1, \eta^3_2$ are semi-basic for the natural 
projection $\calF(\bb{E}^3) \to M$, and the form $\eta^1_2$ spans the 
cotangent space of each fiber of the projection.  While the forms 
$\eta^1, \eta^2, \eta^3_1, \eta^3_2$ are not well-defined on $M$, certain 
combinations of them are.  In particular, since $\eta^3$ is 
well-defined on $M$, so is the form
\[ d\eta^3 = -\eta^3_1 \w \eta^1 - \eta^3_2 \w \ \eta^2. \]
In addition, the area form $\eta^1 \w \eta^2$ is well-defined on $M$, 
as are the 2-forms $\eta^3_1 \w \eta^3_2$ and $\eta^3_1 \w \eta^2 + 
\eta^1 \w \eta^3_2$ which describe Gauss and mean curvature.  So in the 
example above, the ideal $\tilde{\calI}$ on $\calF(\bb{E}^3)$ actually 
projects to a well-defined ideal 
\[ \calI = \{\eta^3,\ d\eta^3, \ \eta^3_1 \w \eta^3_2 + \eta^1 \w \eta^2\} \]
on $M$.  Integral manifolds of this ideal are 
the canonical liftings to $M$ of surfaces in $\bb{E}^3$ with constant 
Gauss curvature $K \equiv -1$.  Furthermore, in this case the pair $(M, 
\calI)$ is a hyperbolic Monge-Amp\`ere system.

These constructions can all be carried out when $\bb{E}^3$ is 
replaced by the space forms $S^3, \bb{H}^3$, by flat Lorentzian space 
(which we will denote $\bb{E}^{2,1}$),
or by Lorentzian space forms of constant sectional curvature 1 
or $-1$ (which we will denote $S^{2,1}$ and $\bb{H}^{2,1}$, respectively). 
In each case the frame bundle 
will be isomorphic to the Lie group of isometries of the underlying 
space form, and the structure equations will vary depending on the group.  
In addition, in the Lorentzian case there will be variations 
depending on whether we are considering spacelike or timelike 
surfaces.  In either case we choose orthonormal frames along the surface 
with $e_1$ and $e_2$ tangent to the surface;
in the spacelike case we choose frames with 
\[ \langle e_1, e_1 \rangle = \langle e_2, e_2 \rangle = 1, \ \ 
\langle e_3, e_3 \rangle = -1, \]
and in the timelike case we choose frames with 
\[ \langle e_1, e_1 \rangle = \langle e_3, e_3 \rangle = 1, \ \ 
\langle e_2, e_2 \rangle = -1. \]
For spacelike surfaces in either Riemannian or Lorentzian space forms, 
the Gauss equation
\[  d\eta^1_2 = K\, \eta^1 \w \eta^2 \]
is taken as a definition of the Gauss curvature of the surface.
For timelike surfaces in Lorentzian space forms, the analog of the Gauss 
equation is 
\[ d\eta^1_2 = -K\, \eta^1 \w \eta^2. \]
(See \cite{O83} for a discussion of curvature in Lorentzian spaces.)
Moreover, whenever the underlying space form
has nonzero sectional curvature $K_0$, the relationship of between 
the Gauss curvature $K$ of a surface and the second fundamental form of the 
surface is
\[ K = K_0 + (h_{11} h_{22} - h_{12}^2) \]
when the underlying space form is Riemannian and
\[ K = K_0 - (h_{11} h_{22} - h_{12}^2) \]
for either spacelike or timelike surfaces when the underlying space 
form is Lorent\-zian.
Thus we have
\[ \eta^3_1 \w \eta^3_2 = (K - K_0)\, \eta^1 \w \eta^2 \]
for surfaces in Riemannian space forms and
\[ \eta^3_1 \w \eta^3_2 = (K_0 - K)\, \eta^1 \w \eta^2 \]
for either spacelike or timelike surfaces in Lorentzian space forms.
Finally, for timelike surfaces in Lorentzian space forms the mean curvature 
is defined to be 
one-half of the trace of $II$ with repect to the Lorentzian metric 
$I$, so 
\[ H = \tfrac{1}{2} (h_{11} - h_{22}). \]
In this case we have
\[ \eta^3_1 \w \eta^2 - \eta^1 \w \eta^3_2 = (h_{11} - h_{22}) \,\eta^1 \w 
\eta^2 = 2H\, \eta^1 \w \eta^2. \]

The structure equations in the various cases are:
\begin{itemize}
\item{Surfaces in $\bb{E}^3$: the frame bundle is isomorphic to 
$E(3)$, and the structure equations are
\begin{align*}
d\eta^i & = - \sum_{j=1}^3 \eta^i_j \w \eta^j \\ 
d\eta^i_j & = - \sum_{k=1}^3 \eta^i_k \w \eta^k_j, \notag
\end{align*}
with $\eta^i_j = -\eta^j_i$.}
\item{Surfaces in $S^3$: the frame bundle is isomorphic to $O(4)$, 
and the structure equations are
\begin{align*}
d\eta^i & = - \sum_{j=1}^3 \eta^i_j \w \eta^j \\ 
d\eta^i_j & = - \sum_{k=1}^3 \eta^i_k \w \eta^k_j + \eta^i \w \eta^j, \notag
\end{align*}
with $\eta^i_j = -\eta^j_i$.}
\item{Surfaces in $\bb{H}^3$: the frame bundle is isomorphic to $O(3,1)$, 
and the structure equations are
\begin{align*}
d\eta^i & = - \sum_{j=1}^3 \eta^i_j \w \eta^j \\ 
d\eta^i_j & = - \sum_{k=1}^3 \eta^i_k \w \eta^k_j - \eta^i \w \eta^j, \notag
\end{align*}
with $\eta^i_j = -\eta^j_i$.}
\item{Spacelike surfaces in $\bb{E}^{2,1}$: the frame bundle is isomorphic to 
the Lorentzian group $E(2,1)$ (i.e., the Lorentzian analog of $E(3)$), 
and the structure equations are
\begin{align*}
d\eta^i & = - \sum_{j=1}^3 \eta^i_j \w \eta^j \\ 
d\eta^i_j & = - \sum_{k=1}^3 \eta^i_k \w \eta^k_j, \notag
\end{align*}
with $\eta^i_i = 0,\ \eta^2_1 = -\eta^1_2, \ \eta^1_3 = \eta^3_1, \ 
\eta^2_3 = \eta^3_2$.}
\item{Timelike surfaces in $\bb{E}^{2,1}$: the frame bundle is isomorphic to 
the Lorentzian group $E(2,1)$, and the structure equations are
\begin{align*}
d\eta^i & = - \sum_{j=1}^3 \eta^i_j \w \eta^j \\ 
d\eta^i_j & = - \sum_{k=1}^3 \eta^i_k \w \eta^k_j, \notag
\end{align*}
with $\eta^i_i = 0,\ \eta^2_1 = \eta^1_2, \ \eta^1_3 = -\eta^3_1, \ 
\eta^2_3 = \eta^3_2$.}
\item{Spacelike surfaces in $S^{2,1}$: the frame 
bundle is isomorphic to $O(3,1)$, and the structure equations are
\begin{align*}
d\eta^i & = - \sum_{j=1}^3 \eta^i_j \w \eta^j \notag \\
d\eta^1_2 & = - \eta^1_3 \w \eta^3_2 + \eta^1 \w \eta^2 
\\ 
d\eta^3_1 & = - \eta^3_2 \w \eta^2_1 + \eta^3 \w \eta^1 \notag \\
d\eta^3_2 & = - \eta^3_1 \w \eta^1_2 + \eta^3 \w \eta^2 \notag
\end{align*}
with $\eta^i_i = 0,\ \eta^2_1 = -\eta^1_2, \ \eta^1_3 = \eta^3_1, \ 
\eta^2_3 = \eta^3_2$.
It is straightforward to show that
this case is actually isomorphic to the case of surfaces in 
$\bb{H}^3$ via the change of basis
\[ \{\eta^1,\, \eta^2,\, \eta^3,\, \eta^1_2,\, \eta^3_1,\, \eta^3_2\}\  \to  
\ \{-\eta^3_2,\, \eta^3_1,\, \eta^3,\, \eta^1_2,\, \eta^2,\, -\eta^1\}. \]
This correspondence sends surfaces of Gauss curvature $K \neq -1$ in 
$\bb{H}^3$ to their Gauss images, which are spacelike 
surfaces of Gauss curvature $K \neq 1$ in $S^{2,1}$.
}
\item{Timelike surfaces in $S^{2,1}$: the frame 
bundle is isomorphic to $O(3,1)$, and the structure equations are
\begin{align*}
d\eta^i & = - \sum_{j=1}^3 \eta^i_j \w \eta^j \notag \\
d\eta^1_2 & = - \eta^1_3 \w \eta^3_2 - \eta^1 \w \eta^2 
\\ 
d\eta^3_1 & = - \eta^3_2 \w \eta^2_1 + \eta^3 \w \eta^1 \notag \\
d\eta^3_2 & = - \eta^3_1 \w \eta^1_2 - \eta^3 \w \eta^2 \notag
\end{align*}
with $\eta^i_i = 0,\ \eta^2_1 = \eta^1_2, \ \eta^1_3 = -\eta^3_1, \ 
\eta^2_3 = \eta^3_2$.}
\item{Spacelike surfaces in $\bb{H}^{2,1}$: the frame 
bundle is isomorphic to $O(2,2)$, and the structure equations are
\begin{align*}
d\eta^i & = - \sum_{j=1}^3 \eta^i_j \w \eta^j \notag \\
d\eta^1_2 & = - \eta^1_3 \w \eta^3_2 - \eta^1 \w \eta^2 
\\ 
d\eta^3_1 & = - \eta^3_2 \w \eta^2_1 - \eta^3 \w \eta^1 \notag \\
d\eta^3_2 & = - \eta^3_1 \w \eta^1_2 - \eta^3 \w \eta^2 \notag
\end{align*}
with $\eta^i_i = 0,\ \eta^2_1 = -\eta^1_2, \ \eta^1_3 = \eta^3_1, \ 
\eta^2_3 = \eta^3_2$.}
\item{Timelike surfaces in $\bb{H}^{2,1}$: the frame 
bundle is isomorphic to $O(2,2)$, and the structure equations are
\begin{align*}
d\eta^i & = - \sum_{j=1}^3 \eta^i_j \w \eta^j \notag \\
d\eta^1_2 & = - \eta^1_3 \w \eta^3_2 + \eta^1 \w \eta^2 
\\ 
d\eta^3_1 & = - \eta^3_2 \w \eta^2_1 - \eta^3 \w \eta^1 \notag \\
d\eta^3_2 & = - \eta^3_1 \w \eta^1_2 + \eta^3 \w \eta^2 \notag
\end{align*}
with $\eta^i_i = 0,\ \eta^2_1 = \eta^1_2, \ \eta^1_3 = -\eta^3_1, \ 
\eta^2_3 = \eta^3_2$.}
\end{itemize}

\section{Interpretation of Cases 3C and 3D}

In Cases 3C and 3D, we found a coframing $\{\theta_1, \theta_2, 
\om^1, \om^2, \om^3, \om^4\}$ whose structure equations have constant 
coefficients.  This implies that the forms in the coframing form a Lie 
algebra.  This in turn gives the manifold $\calB$ a Lie group 
structure (at least 
locally) by regarding the forms in the coframing as the 
left-invariant forms on $\calB$.  The first step in 
interpreting the structure equations is to identify the Lie algebra 
that they define, and in all but one case it turns out 
to be one of those described in the previous section.  Then because 
the contact forms $\theta_1, \theta_2$ are each determined up to 
scalar multiples, we must find two distinct bases for the Lie 
algebra: a basis $\{\eta^i, \eta^i_j\}$ for which $\eta^3$ is a multiple 
of $\theta_1$, and a basis  $\{\zeta^i, \zeta^i_j\}$
for which $\zeta^3$ is a multiple of $\theta_2$.  The B\"acklund 
transformation is then given by the transformation relating these two 
bases for the Lie algebra.  These transformations can all be described 
by geodesic congruences of some sort, in the same way that the 
classical B\"acklund transformation between pseudospherical surfaces 
is given by line congruences.

These computations were carried out using the algorithm in \cite{RWZ88} 
with the assistance of Maple.  The algorithm divides into several 
cases depending on the value of $B_1$ in case 3C and the values of 
$B_1, B_3, \eps$ in case 3D.  The change-of-basis matrices are rather 
complicated and not very enlightening, so they will be omitted here. 

\subsection{Case 3C}

Recall that the structure equations in this case are 
\begin{align*}
d\theta_1 & = \theta_1 \w (\om^1 + \om^3) +  \om^1 \w \om^2 + \om^3 \w 
\om^4 \\
d\theta_2 & = -\theta_2 \w (\om^1 + \om^3) +  \om^1 \w \om^2 - \om^3 \w 
\om^4 \\
d\om^1 & = B_1 (\theta_1 \w \theta_2 + \theta_1 \w \om^2 + \theta_2 
\w \om^2 + \om^2 \w \om^3) \\
d\om^2 & = \frac{1}{B_1}\, (\theta_1 + \theta_2 - \om^3) \w \om^1 + 
\om^3 \w \om^4 \\
d\om^3 & = (\theta_1 - \theta_2 - B_1\, \om^2) \w \om^3 \\
d\om^4 & = \theta_1 \w \theta_2 - \theta_1 \w \om^4 + \theta_2 \w 
\om^4 + B_1\, \om^2 \w \om^4 + \om^1 \w \om^2
\end{align*}
with $B_1 \neq 0$.
Carrying out the algorithm described above shows that:
\begin{itemize}
\item{If $B_1 \neq 
2$, then the Lie algebra is $\mathfrak{so}(2,2)$.  For each of the 
two bases computed by the algorithm, the structure 
equations coincide with those for timelike surfaces in $\bb{H}^{2,1}$; 
moreover, the ideals $\calI_1, \calI_2$ take the form
\begin{gather*}
\calI_1 = \{\eta^3,\ d\eta^3,\ \eta^3_1 \w \eta^2 - \eta^1 \w 
\eta^3_2 - 2 \eta^1 \w \eta^2 \} \\
\calI_2 = \{\zeta^3,\ d\zeta^3,\ \zeta^3_1 \w \zeta^2 - \zeta^1 \w 
\zeta^3_2 - 2 \zeta^1 \w \zeta^2 \}.
\end{gather*}
So up to contact equivalence, $\calB$ represents a transformation between 
timelike surfaces of 
constant mean curvature equal to 1 in $\bb{H}^{2,1}$.  We note that the 
change-of-basis matrices have different 
expressions for $B_1$ in the ranges $B_1 < 0,\ \ 0 < B_1 < 2$, and  $B_1 > 
2$.
}
\item{If $B_1 = 2$, then the Lie algebra is $\mathfrak{e}(2,1)$.
For each of the two bases computed by the algorithm, the structure 
equations coincide with those for timelike surfaces in $\bb{E}^{2,1}$; 
moreover, the ideals $\calI_1, \calI_2$ take the form
\begin{gather*}
\calI_1 = \{\eta^3, \ d\eta^3, \ \eta^3_1 \w \eta^2 - \eta^1 \w 
\eta^3_2 \} \\
\calI_2 = \{\zeta^3,\ d\zeta^3,\ \zeta^3_1 \w \zeta^2 - \zeta^1 \w 
\zeta^3_2 \}.
\end{gather*}
So up to contact equivalence, $\calB$ represents a transformation between 
timelike minimal surfaces in $\bb{E}^{2,1}$.  This transformation is 
explored in detail in \cite{C01}.}
\end{itemize}

Thus we have the following theorem.

\begin{theorem}\label{case3cthm}
Let $\calB \subset M_1 \times M_2$ be a homogeneous B\"acklund 
transformation with the vectors $[C_1 \ \ C_2], \ [C_3 \ \ C_4], \ 
[B_1 \ \ B_2], \ [B_3 \ \ B_4]$ all nonzero, the pair $[C_1 \ \ 
C_2], \ [B_1 \ \ B_2]$ linearly independent, and the pair $[C_3 \ \ C_4], 
\ [B_3 \ \ B_4]$ linearly dependent.  Then $\calB$ is locally contact 
equivalent to either
\begin{enumerate}
\item{A B\"acklund transformation between timelike minimal surfaces in 
$\bb{E}^{2,1}$, or}
\item{A B\"acklund transformation between timelike surfaces of 
constant mean curvature equal to 1 in $\bb{H}^{2,1}$.}
\end{enumerate}
In both cases, the transformation may be described in terms of geodesic 
congruences.
\end{theorem}

\subsection{Case 3D}

Recall that the structure equations in this case are
\begin{align*}
d\theta_1 & = \theta_1 \w (\om^1 + \om^3) +  
\om^1 \w \om^2 + \om^3 \w \om^4 \\
d\theta_2 & = -\theta_2 \w (\om^1 + \om^3) +  
\om^1 \w \om^2 - \om^3 \w \om^4 \\
d\om^1 & = (B_1\, \theta_1 + B_1\, \theta_2 - 
\eps\, \om^4) \w \om^2 + B_1\, \theta_1 \w \theta_2 \\
d\om^2 & = (-\eps B_3\, \theta_1 - \eps B_3\, 
\theta_2 + \frac{B_3}{B_1}\, \om^4) \w \om^1 + \om^3 \w \om^4  \\
d\om^3 & = (-B_3\, \theta_1 + B_3\, \theta_2 - \eps\, \om^2) 
\w \om^4 + B_3\, \theta_1 \w \theta_2  \\
d\om^4 & = (\eps B_1\, \theta_1 - \eps B_1\, 
\theta_2 + \frac{B_1}{B_3}\, \om^2) \w \om^3 + \om^1 \w \om^2
\end{align*}
with $B_1, B_3 \neq 0$ and $\eps = \pm 1$.  The algorithm described 
above divides into many cases depending on the values of these 
parameters.

When $\eps = -1$, the $B_1 B_3$ plane divides into regions as shown 
in Figure 1.  The curve in this graph is defined by the equation
\[ 4 B_1^2 B_3^2 - 4 B_1^2 B_3 + 4 B_1 B_3^2 + B_1^2 + 2 B_1 B_3 + 
B_3^2 = 0, \]
and it may be parametrized by
\[ B_1 = -\frac{1}{2} (t + 1)^2, \qquad B_3 = \frac{1}{2} 
\left(\frac{1}{t} + 1\right)^2 \]
for $t \neq 0$.  (The point corresponding to $t=-1$ is (0,0) and so is 
not included in our parameter space.)  For convenience, we define
\[ Q^- = 4 B_1^2 B_3^2 - 4 B_1^2 B_3 + 4 B_1 B_3^2 + B_1^2 + 2 B_1 B_3 + 
B_3^2. \]
When the point $(B_1, B_3)$ is in the second or fourth quadrant, 
$Q^-$ can be factored as
\[ Q^- = q^-_1 q^-_2 \]
with
\begin{align*}
q^-_1 &= 2 B_1 B_3 - B_1 + B_3 + 2\sqrt{-B_1 B_3} \\
q^-_2 &= 2 B_1 B_3 - B_1 + B_3 - 2\sqrt{-B_1 B_3}.
\end{align*}
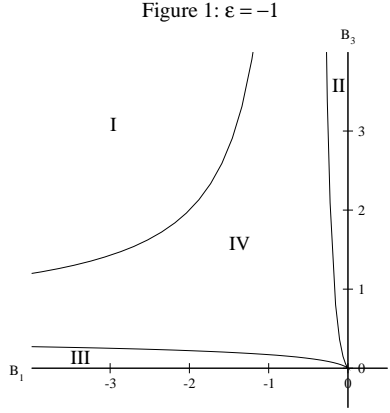
\begin{center}
\epsfig{figure=homogfig1.eps,width=.4\linewidth}
\end{center}

\begin{itemize}
\item{If $Q^- = 0$ and $t > 0$ (so that $B_1 < -\tfrac{1}{2}$ and $B_3 > 
\tfrac{1}{2}$), then the Lie algebra is $\mathfrak{e}(3)$.
For each of the two bases computed by the algorithm, the structure 
equations coincide with those for surfaces in $\bb{E}^3$; 
moreover, the ideals $\calI_1, \calI_2$ take the form
\begin{gather*}
\calI_1 = \{\eta^3,\ d\eta^3,\ \eta^3_1 \w \eta^3_2  + \eta^1 \w 
\eta^2 \} \\
\calI_2 = \{\zeta^3,\ d\zeta^3,\ \zeta^3_1 \w \zeta^3_2 + \zeta^1 \w 
\zeta^2 \}.
\end{gather*}
So up to contact equivalence, $\calB$ represents a transformation 
between surfaces of constant Gauss curvature $K = -1$ in $\bb{E}^3$.  This 
is the classical B\"acklund transformation between 
pseudospherical surfaces, and the parameter $t$ along the curve $Q^- = 
0$ is a function of the usual parameter appearing in this 
transformation.}
\item{If $Q^-=0$ and $t < 0$ (so that either $B_1 > -\tfrac{1}{2}$ or 
$B_3 < \tfrac{1}{2}$), then the Lie algebra is $\mathfrak{e}(2,1)$.
For each of the two bases computed by the algorithm, the structure 
equations coincide with those for spacelike surfaces in $\bb{E}^{2,1}$; 
moreover, the ideals $\calI_1, \calI_2$ take the form
\begin{gather*}
\calI_1 = \{\eta^3,\ d\eta^3,\ \eta^3_1 \w \eta^3_2  + \eta^1 \w 
\eta^2 \} \\
\calI_2 = \{\zeta^3,\ d\zeta^3,\ \zeta^3_1 \w \zeta^3_2 + \zeta^1 \w 
\zeta^2 \}.
\end{gather*}
So up to contact equivalence, $\calB$ represents a transformation 
between spacelike surfaces of constant Gauss curvature $K = 1$ in 
$\bb{E}^{2,1}$.  
}
\item{In Region I of Quadrant 2, the Lie algebra is 
$\mathfrak{so}(4)$.
For each of the two bases computed by the algorithm, the structure 
equations coincide with those for surfaces in $S^3$; 
moreover, the ideals $\calI_1, \calI_2$ take the form
\begin{gather*}
\calI_1 = \{\eta^3,\ d\eta^3,\ \eta^3_1 \w \eta^3_2  + 
\left(\frac{q^-_1}{q^-_2}\right) \eta^1 \w \eta^2 \} \\
\calI_2 = \{\zeta^3,\ d\zeta^3,\ \zeta^3_1 \w \zeta^3_2 + 
\left(\frac{q^-_1}{q^-_2}\right) \zeta^1 \w \zeta^2 \}.
\end{gather*}
So up to contact equivalence, $\calB$ represents a transformation 
between surfaces of constant Gauss curvature
\[ K = 1 -\frac{q^-_1}{q^-_2} \] 
in $S^3$.  As $(B_1, B_3)$ ranges over 
Region I, $K$ takes values in the interval $(0, 1)$.}
\item{In Regions II and III of Quadrant 2 and in Quadrant 4, the Lie 
algebra is $\mathfrak{so}(2,2)$.
For each of the two bases computed by the algorithm, the structure 
equations coincide with those for spacelike surfaces in $\bb{H}^{2,1}$; 
moreover, the ideals $\calI_1, \calI_2$ take the form
\begin{gather*}
\calI_1 = \{\eta^3,\ d\eta^3,\ \eta^3_1 \w \eta^3_2  + 
\left(\frac{q^-_2}{q^-_1}\right) \eta^1 \w \eta^2 \} \\
\calI_2 = \{\zeta^3,\ d\zeta^3,\ \zeta^3_1 \w \zeta^3_2 + 
\left(\frac{q^-_2}{q^-_1}\right) \zeta^1 \w \zeta^2 \}.
\end{gather*}
So up to contact equivalence, $\calB$ represents a transformation 
between spacelike surfaces of constant Gauss curvature
\[ K = \frac{q^-_2}{q^-_1} - 1\] 
in $\bb{H}^{2,1}$.  As $(B_1, B_3)$ ranges over 
these regions, $K$ takes values in the interval $(-1, 0)$ in Regions 
II and III of Quadrant 2 and in the interval $(0, \infty)$ 
in Quadrant 4.
We note that the change-of-basis matrices have 
different expressions in each of the three regions.}
\item{In Region IV of Quadrant 2, the Lie algebra is 
$\mathfrak{so}(3,1)$.
For each of the two bases computed by the algorithm, the structure 
equations coincide with those for surfaces in $\bb{H}^3$, or 
equivalently, for spacelike surfaces in $S^{2,1}$.  Regarded as surfaces 
in $\bb{H}^3$, the ideals $\calI_1, \calI_2$ take the form
\begin{gather*}
\calI_1 = \{\eta^3,\ d\eta^3,\ \eta^3_1 \w \eta^3_2  -
\left(\frac{q^-_2}{q^-_1}\right) \eta^1 \w \eta^2 \} \\
\calI_2 = \{\zeta^3,\ d\zeta^3,\ \zeta^3_1 \w \zeta^3_2 - 
\left(\frac{q^-_2}{q^-_1}\right) \zeta^1 \w \zeta^2 \}.
\end{gather*}
Regarded as spacelike surfaces in $S^{2,1}$, the ideals $\calI_1, 
\calI_2$ take the form
\begin{gather*}
\calI_1 = \{\eta^3,\ d\eta^3,\ \eta^3_1 \w \eta^3_2  - 
\left(\frac{q^-_1}{q^-_2}\right) \eta^1 \w \eta^2 \} \\
\calI_2 = \{\zeta^3,\ d\zeta^3,\ \zeta^3_1 \w \zeta^3_2 -
\left(\frac{q^-_1}{q^-_2}\right) \zeta^1 \w \zeta^2 \}.
\end{gather*}
So up to contact equivalence, $\calB$ may be regarded as representing 
either a transformation between surfaces of constant Gauss curvature
\[ K = \frac{q^-_2}{q^-_1} - 1\] 
in $\bb{H}^3$, or a transformation between spacelike surfaces of 
constant Gauss curvature
\[ K = 1 - \frac{q^-_1}{q^-_2} \]
in $S^{2,1}$.  In the first case $K$ takes values in the interval 
$(-\infty, -1)$ as $(B_1, B_3)$ ranges over Region IV, and in the 
second case $K$ takes values in the interval $(1, \infty)$ as $(B_1, 
B_3)$ ranges over Region IV.}
\item{In Quadrants 1 and 3, the Lie algebra is $\mathfrak{so}(2,2)$.
For each of the two bases computed by the algorithm, the structure 
equations coincide with those for timelike surfaces in $\bb{H}^{2,1}$; 
moreover, the ideals $\calI_1, \calI_2$ take the form
\begin{gather*}
\calI_1 = \{\eta^3,\ d\eta^3,\ \eta^3_1 \w \eta^2 - \eta^1 \w 
\eta^3_2 - 2\frac{2 B_1 B_3 - B_1 + B_3}{\sqrt{Q^-}} \eta^1 \w \eta^2 \} \\
\calI_2 =\{\zeta^3,\ d\zeta^3,\ \zeta^3_1 \w \zeta^2 - \zeta^1 \w 
\zeta^3_2 - 2\frac{2 B_1 B_3 - B_1 + B_3}{\sqrt{Q^-}} \zeta^1 \w \zeta^2 \}. 
\end{gather*}
So up to contact equivalence, $\calB$ represents a transformation 
between timelike surfaces of constant mean curvature
\[ H = \frac{2 B_1 B_3 - B_1 + B_3}{\sqrt{Q^-}}\]
in $\bb{H}^{2,1}$.  As $(B_1, B_3)$ ranges over 
these regions, $H$ takes values in the interval $(-1, 1)$.
We note that the change-of-basis matrices have 
different expressions in each quadrant.}
\end{itemize}

When $\eps = 1$, the $B_1 B_3$ plane divides into regions as shown 
in Figure 2.  The curve in this graph is defined by the equation
\[ 4 B_1^2 B_3^2 + 4 B_1^2 B_3 - 4 B_1 B_3^2 + B_1^2 + 2 B_1 B_3 + 
B_3^2 = 0, \]
and it may be parametrized by
\[ B_1 = \frac{1}{2} (t + 1)^2, \qquad B_3 = -\frac{1}{2} 
\left(\frac{1}{t} + 1\right)^2 \]
for $t \neq 0$.  (The point corresponding to $t=-1$ is (0,0) and so is 
not included in our parameter space.)  For convenience, we define
\[ Q^+ = 4 B_1^2 B_3^2 + 4 B_1^2 B_3 - 4 B_1 B_3^2 + B_1^2 + 2 B_1 B_3 + 
B_3^2. \]
When the point $(B_1, B_3)$ is in the second or fourth quadrant, 
$Q^+$ can be factored as
\[ Q^+ = q^+_1 q^+_2 \]
with
\begin{align*}
q^+_1 &= 2 B_1 B_3 + B_1 - B_3 + 2\sqrt{-B_1 B_3} \\
q^+_2 &= 2 B_1 B_3 + B_1 - B_3 - 2\sqrt{-B_1 B_3}.
\end{align*}
\begin{center}
\epsfig{figure=homogfig2.eps,width=.4\linewidth}
\end{center}

\begin{itemize}
\item{If $Q^+ = 0$ then the Lie algebra is $\mathfrak{e}(2,1)$.
For each of the two bases computed by the algorithm, the structure 
equations coincide with those for timelike surfaces in $\bb{E}^{2,1}$; 
moreover, the ideals $\calI_1, \calI_2$ can be written either in the form
\begin{gather*}
\calI_1 = \{\eta^3,\ d\eta^3,\ \eta^3_1 \w \eta^3_2  + \eta^1 \w 
\eta^2 \} \\
\calI_2 = \{\zeta^3,\ d\zeta^3,\ \zeta^3_1 \w \zeta^3_2 + \zeta^1 \w 
\zeta^2 \}
\end{gather*}
or in the form
\begin{gather*}
\calI_1 = \{\eta^3,\ d\eta^3,\ \eta^3_1 \w \eta^2 - \eta^1 \w \eta^3_2  - 
2\eta^1 \w \eta^2 \} \\
\calI_2 = \{\zeta^3,\ d\zeta^3,\ \zeta^3_1 \w \zeta^2 - \zeta^1 \w 
\zeta^3_2 - 2 \zeta^1 \w \zeta^2 \} .
\end{gather*}
So up to contact equivalence, $\calB$ may be regarded as representing 
either a transformation between timelike surfaces of constant Gauss 
curvature $K = 1$ or a transformation between timelike surfaces of 
constant mean curvature $H = 1$ in $\bb{E}^{2,1}$.}
\item{In Region I of Quadrant 4, the Lie 
algebra is $\mathfrak{so}(2,2)$.
For each of the two bases computed by the algorithm, the structure 
equations coincide with those for timelike surfaces in $\bb{H}^{2,1}$; 
moreover, the ideals $\calI_1, \calI_2$ can be written either in the form
\begin{gather*}
\calI_1 = \{\eta^3,\ d\eta^3,\ \eta^3_1 \w \eta^3_2  + 
\left(\frac{q^+_1}{q^+_2}\right) \eta^1 \w \eta^2 \} \\
\calI_2 = \{\zeta^3,\ d\zeta^3,\ \zeta^3_1 \w \zeta^3_2 + 
\left(\frac{q^+_1}{q^+_2}\right) \zeta^1 \w \zeta^2 \}
\end{gather*}
or in the form
\begin{gather*}
\calI_1 = \{\eta^3,\ d\eta^3,\ \eta^3_1 \w \eta^2 - \eta^1 \w \eta^3_2  - 
2 \left( \frac{2 B_1 B_3 + B_1 - B_3}{\sqrt{Q^+}} \right)\,  \eta^1 \w 
\eta^2 \}  \\
\calI_2 = \{\zeta^3,\ d\zeta^3,\ \zeta^3_1 \w \zeta^2 - \zeta^1 \w 
\zeta^3_2 - 2 \left( \frac{2 B_1 B_3 + B_1 - B_3}{\sqrt{Q^+}}\right) \, 
\zeta^1 \w \zeta^2 \} .
\end{gather*}
So up to contact equivalence, $\calB$ may be regarded as representing 
either a transformation between timelike surfaces of constant Gauss 
curvature 
\[ K =  \frac{q^+_1}{q^+_2} - 1 \] 
or a transformation between timelike surfaces of 
constant mean curvature 
\[ H = \frac{2 B_1 B_3 + B_1 - B_3}{\sqrt{Q^+}} \]
in $\bb{H}^{2,1}$.  As $(B_1, B_3)$ ranges over 
Region I, $K$ takes values in the interval $(0, 1)$; meanwhile, $H$ takes 
values in the interval $(-\infty, -1)$.}
\item{In Regions II and III of Quadrant 4 and in Quadrant 2, the Lie 
algebra is $\mathfrak{so}(2,2)$.
For each of the two bases computed by the algorithm, the structure 
equations coincide with those for timelike surfaces in $\bb{H}^{2,1}$; 
moreover, the ideals $\calI_1, \calI_2$ can be written either in the form
\begin{gather*}
\calI_1 = \{\eta^3,\ d\eta^3,\ \eta^3_1 \w \eta^3_2  + 
\left(\frac{q^+_2}{q^+_1}\right)  \eta^1 \w \eta^2 \} \\
\calI_2 = \{\zeta^3,\ d\zeta^3,\ \zeta^3_1 \w \zeta^3_2 + 
\left(\frac{q^+_2}{q^+_1}\right)  \zeta^1 \w \zeta^2 \}
\end{gather*}
or in the form
\begin{gather*}
\calI_1 = \{\eta^3,\ d\eta^3,\ \eta^3_1 \w \eta^2 - \eta^1 \w \eta^3_2  - 
2 \left( \frac{2 B_1 B_3 + B_1 - B_3}{\sqrt{Q^+}}\right) \,  \eta^1 \w 
\eta^2 \} \\
\calI_2 = \{\zeta^3,\ d\zeta^3,\ \zeta^3_1 \w \zeta^2 - \zeta^1 \w 
\zeta^3_2 - 2 \left( \frac{2 B_1 B_3 + B_1 - B_3}{\sqrt{Q^+}}\right) \, 
\zeta^1 \w \zeta^2 \} .
\end{gather*}
So up to contact equivalence, $\calB$ may be regarded as representing 
either a transformation between timelike surfaces of constant Gauss 
curvature 
\[ K =  \frac{q^+_2}{q^+_1} - 1\] 
or a transformation between timelike surfaces of 
constant mean curvature 
\[ H = \frac{2 B_1 B_3 + B_1 - B_3}{\sqrt{Q^+}} \]
in $\bb{H}^{2,1}$.
As $(B_1, B_3)$ ranges over 
these regions, $K$ takes values in the interval $(-1, 0)$ in Regions 
II and III of Quadrant 4 and in the interval $(0, \infty)$ 
in Quadrant 2; meanwhile, $H$ takes values in the interval $(1, 
\infty)$ in Regions II and III of Quadrant 4 and in the interval $(-\infty, 
-1)$ in Quadrant 2.  We note that the change-of-basis matrices have 
different expressions in each of the three regions.}  
\item{In Region IV of Quadrant 4, the Lie algebra is 
$\mathfrak{so}(3,1)$.
For each of the two bases computed by the algorithm, the structure 
equations coincide with those for timelike surfaces in $S^{2,1}$; 
moreover, the ideals $\calI_1, \calI_2$ can be written either in the form
\begin{gather*}
\calI_1 = \{\eta^3,\ d\eta^3,\ \eta^3_1 \w \eta^3_2  - 
\left(\frac{q^+_2}{q^+_1}\right)  \eta^1 \w \eta^2 \} \\
\calI_2 = \{\zeta^3,\ d\zeta^3,\ \zeta^3_1 \w \zeta^3_2 - 
\left(\frac{q^+_2}{q^+_1}\right)  \zeta^1 \w \zeta^2 \}
\end{gather*}
or in the form
\begin{gather*}
\calI_1 = \{\eta^3,\ d\eta^3,\ \eta^3_1 \w \eta^2 - \eta^1 \w \eta^3_2  - 
2 \left( \frac{2 B_1 B_3 + B_1 - B_3}{\sqrt{-Q^+}}\right) \,  \eta^1 \w 
\eta^2  \} \\
\calI_2 = \{\zeta^3,\ d\zeta^3,\ \zeta^3_1 \w \zeta^2 - \zeta^1 \w 
\zeta^3_2 - 2 \left( \frac{2 B_1 B_3 + B_1 - B_3}{\sqrt{-Q^+}}\right) \, 
\zeta^1 \w \zeta^2 \} .
\end{gather*}
So up to contact equivalence, $\calB$ may be regarded as representing 
either a transformation between timelike surfaces of constant Gauss 
curvature 
\[ K = 1 - \frac{q^+_2}{q^+_1}\] 
or a transformation between timelike surfaces of 
constant mean curvature 
\[ H = \frac{2 B_1 B_3 + B_1 - B_3}{\sqrt{-Q^+}} \]
in $S^{2,1}$.  As $(B_1, B_3)$ ranges over 
Region IV, $K$ takes values in the interval $(1, \infty)$; meanwhile, 
$H$ ranges over all real numbers.}  
\item{In Quadrants 1 and 3, the Lie algebra is $\mathfrak{so}(3) 
\oplus \mathfrak{so}(2,1)$.  The corresponding Lie group is denoted 
$SO^{\ast}(4)$ in Cartan's list of Lie groups as described in 
\cite{H62}.  This group has no 
natural 3-dimensional quotients compatible with the contact 
structures given by $\theta_1$ and $\theta_2$, and so there is no natural 
way to regard these examples as transformations of surfaces in any
3-dimensional space.  They may naturally be regarded as transformations of 
certain surfaces in a 5-dimensional quotient space of $SO^{\ast}(4)$.}
\end{itemize}

Putting all these cases together yields the following theorem.

\begin{theorem}\label{case3dthm}
Let $\calB \subset M_1 \times M_2$ be a homogeneous B\"acklund 
transformation with the vectors $[C_1 \ \ C_2], \ [C_3 \ \ C_4], \ 
[B_1 \ \ B_2], \ [B_3 \ \ B_4]$ all nonzero and the pairs $[C_1 \ \ 
C_2], \ [B_1 \ \ B_2]$ and $[C_3 \ \ C_4], 
\ [B_3 \ \ B_4]$ both linearly independent.  Then $\calB$ is locally contact 
equivalent to one of the following:
\begin{enumerate}
\item{A B\"acklund transformation between surfaces of constant 
negative Gauss curvature in $\bb{E}^3$}
\item{A B\"acklund transformation between surfaces of constant 
Gauss curvature $0 < K < 1$ in $S^3$}
\item{A B\"acklund transformation between surfaces of constant 
Gauss curvature $-\infty < K < -1$ in $\bb{H}^3$}
\item{A B\"acklund transformation between spacelike surfaces of 
constant positive Gauss curvature in $\bb{E}^{2,1}$}
\item{A B\"acklund transformation between timelike surfaces of 
constant positive \linebreak Gauss curvature, or equivalently, constant nonzero 
mean curvature, in $\bb{E}^{2,1}$}
\item{A B\"acklund transformation between spacelike surfaces of 
constant Gauss curvature $1 < K < \infty$ in $S^{2,1}$}
\item{A B\"acklund transformation between timelike surfaces of 
constant Gauss curvature $1 < K < \infty$, or equivalently, constant 
mean curvature $H \in \bb{R}$, in $S^{2,1}$}
\item{A B\"acklund transformation between spacelike surfaces of 
constant Gauss curvature $-1 < K < \infty, \ K \neq 0$ in $\bb{H}^{2,1}$}
\item{A B\"acklund transformation between timelike surfaces of 
constant Gauss curvature $-1 < K < \infty, \ K \neq 0$, or equivalently, 
constant mean curvature $|H| > 1$, in $\bb{H}^{2,1}$}
\item{A B\"acklund transformation between timelike surfaces of 
constant mean curvature $|H| < 1$ in $\bb{H}^{2,1}$.}
\item{A B\"acklund transformation between certain surfaces in a 
5-dimensional quotient space of $SO^{\ast}(4)$.}
\end{enumerate}
In all cases, the transformation may be described in terms of geodesic 
congruences.
\end{theorem}

\section{Conclusion}

Theorems \ref{case1thm}, \ref{case2thm}, \ref{case3athm}, \ref{case3bthm}, 
\ref{case3cthm}, and \ref{case3dthm} may be combined to yield the 
following result.

\begin{theorem}\label{bigthm}
Let $\calB \subset M_1 \times M_2$ be a homogeneous B\"acklund 
transformation.  Then $\calB$ is locally contact 
equivalent to one of the following:
\begin{enumerate}
\item{A B\"acklund transformation between solutions of the wave 
equation $z_{xy} = 0$}
\item{A holonomic B\"acklund transformation of the form described in Theorem 
\ref{case3athm}}
\item{The classical B\"acklund transformation between the wave 
equation $z_{xy} = 0$ and Liouville's equation $z_{xy} = e^z$}
\item{A B\"acklund transformation between surfaces of constant 
negative Gauss curvature in $\bb{E}^3$}
\item{A B\"acklund transformation between surfaces of constant 
Gauss curvature $0 < K < 1$ in $S^3$}
\item{A B\"acklund transformation between surfaces of constant 
Gauss curvature $-\infty < K < -1$ in $\bb{H}^3$}
\item{A B\"acklund transformation between spacelike surfaces of 
constant positive Gauss curvature in $\bb{E}^{2,1}$}
\item{A B\"acklund transformation between timelike surfaces of 
constant positive \linebreak Gauss curvature, or equivalently, constant nonzero 
mean curvature, in $\bb{E}^{2,1}$}
\item{A B\"acklund transformation between timelike minimal surfaces 
in $\bb{E}^{2,1}$}
\item{A B\"acklund transformation between spacelike surfaces of 
constant Gauss curvature $1 < K < \infty$ in $S^{2,1}$}
\item{A B\"acklund transformation between timelike surfaces of 
constant Gauss curvature $1 < K < \infty$, or equivalently, constant 
mean curvature $H \in \bb{R}$, in $S^{2,1}$}
\item{A B\"acklund transformation between spacelike surfaces of 
constant Gauss curvature $-1 < K < \infty, \ K \neq 0$ in $\bb{H}^{2,1}$}
\item{A B\"acklund transformation between timelike surfaces of 
constant Gauss curvature $-1 < K < \infty, \ K \neq 0$, or equivalently, 
constant mean curvature $|H| > 1$, in $\bb{H}^{2,1}$}
\item{A B\"acklund transformation between timelike surfaces of 
constant mean curvature $|H| \leq 1$ in $\bb{H}^{2,1}$.}
\item{A B\"acklund transformation between certain surfaces in a 
5-dimensional quotient space of $SO^{\ast}(4)$.}
\end{enumerate}
\end{theorem}

Now this is certainly not the end of the story.  There are 
interesting B\"acklund transformations which are not homogeneous; in 
particular, the classical B\"acklund transformation 
for the sine-Gordon equation
does not appear on this list.  Moreover, the notion of B\"acklund 
transformation used here does not take into account the presence of 
the arbitrary parameter $\lambda$ that plays such an important role 
in the theory of B\"acklund transformations of integrable systems such 
as the sine-Gordon equation.  We hope to address these and other issues in 
future papers.

\end{document}